\documentclass[journal]{IEEEtran}
\usepackage[pdftex]{graphicx}
\usepackage{amsmath,amssymb}
\usepackage{color,cite}
\usepackage{graphicx}
\usepackage{booktabs}
\usepackage{enumitem}
\usepackage{amsfonts}
\usepackage{stmaryrd}
\usepackage{cancel}
\usepackage{mathtools}
\usepackage[hyphens]{url}
\usepackage{subcaption}[font=small]
\usepackage[font=small]{caption}
\usepackage{float}
\usepackage[hidelinks]{hyperref}
\usepackage{cleveref}

\begin{document}

\title{A Hierarchical Local Electricity Market for a DER-rich Grid Edge}

\author{Vineet~Jagadeesan~Nair,~\IEEEmembership{Student Member,~IEEE,} Venkatesh~Venkataramanan,~\IEEEmembership{Member,~IEEE,} Rabab~Haider,~\IEEEmembership{Graduate Student Member,~IEEE,} and 
        Anuradha M. Annaswamy~\IEEEmembership{Fellow,~IEEE}
\thanks{This work was supported by the Department of Energy under Award Number DOE-OE0000920, and Award Number DE-IA0000025. The first author gratefully acknowledges Siemens Corporation for supporting him with a summer internship. The views and opinions of authors expressed herein do not necessarily state or reflect those of the United States Government or any agency thereof. V.~Nair, R.~Haider, and A.M.~Annaswamy are with the Department
of Mechanical Engineering, Massachusetts Institute of Technology (MIT), Cambridge, USA. V.~Venkataramanan is with the National Renewable Energy Laboratory (NREL), e-mail: \{jvineet9, rhaider, aanna\}@mit.edu, venkatesh.venkataramanan@nrel.gov.}} 

\maketitle

\begin{abstract}
With increasing penetration of distributed energy resources (DER) in the distribution system, it is critical to design market structures that enable smooth integration of DERs.  A hierarchical local electricity market (LEM) structure is proposed in this paper with a secondary market (SM) at the lower level representing secondary feeders and a primary market (PM) at the upper level, representing primary feeders, in order to effectively use DERs to increase grid efficiency and resilience. The lower level SM enforces budget, power balance, and flexibility constraints and accounts for costs related to consumers, such as their disutility, flexibility limits, and commitment reliability, while the upper level PM enforces grid physics constraints such as power balance and capacity limits, and also minimizes line losses. The hierarchical LEM is extensively evaluated using a modified IEEE-123 bus with high DER penetration, with each primary feeder consisting of at least three secondary feeders. Data from a GridLAB-D model is used to emulate realistic power injections and load profiles over the course of 24 hours. The performance of the LEM is illustrated by delineating the family of power-injection profiles across the primary and secondary feeders as well as corresponding local electricity tariffs that vary across the distribution grid. Through numerical simulations, the hierarchical LEM is shown to improve the efficiency of the market in terms of lowering overall costs, including both the distribution-level locational marginal prices (d-LMP) as well as retail tariffs paid by customers. Together, it represents an overall framework for a Distribution System Operator (DSO) who can provide the oversight for the entire LEM.
\end{abstract}


\section{Introduction \label{sec:intro}}

The current electrical grid together with the electricity market, was designed for unidirectional energy flows — from large centralized producers to final consumers. However, as the power grid becomes increasingly decentralized, the electricity market structure needs to (1) incentivize consumers appropriately for investing in DERs, (2) ensure that DERs are appropriately compensated for any grid services provided, (3) coordinate with the existing market structure at the wholesale level, (4) provide a suitable rate structure to ensure that relevant costs to various stakeholders are minimized and (5) be scalable with increasing penetration of DERs. Local electricity markets (LEMs) have the potential to empower the consumer to take control of their energy footprint, allow transactive energy trading among members of a community, improve community resilience against wider grid events, and potentially reduce energy bills. Many of the existing consumer choice programs compensate DERs for their grid services, typically through direct incentives and feed-in tariffs. However, these policies do not price the fine-grain locational and temporal variation in the services that DERs are capable of providing, and are therefore unable to meet network-level objectives under high DER penetration \cite{adapen}. As the number of DERs and prosumers (energy consumers who also produce some electricity) increases, more structure is warranted to coordinate them. New operational entities such as a distribution system operator (DSO) that oversees the LEM \cite{renani_optimal_2018}, and new types of tariffs such as distribution-level locational marginal prices (d-LMP), may be needed \cite{wang_distribution_2019}.

Local energy markets have the capability to allow electricity prices to be endogenous quantities rather than be imposed exogenously. In such a marketplace, prosumers can buy and sell energy in an open marketplace, or through an operator \cite{review_LEM_CEM}. The introduction of DERs such as rooftop solar panels and electric vehicles (EV) has introduced significant complexity to the management of the grid. Grid operators and utilities often rely on standard load profiles derived from historical data to model home energy usage and estimate the amount of energy required to supply and balance the grid. However, the intermittent and highly variable nature of the generation from photovoltaic (PV) panels, demand of EVs, and needs of other DERs can cause unpredictable swings in demand. LEMs have the potential to help solve this problem for energy retailers and other grid management entities by offering flexibility services and the opportunity for new business models. LEMs also provide an attractive alternative to FERC Order 2222, as the direct participation of DERs at the wholesale level may introduce tier-bypassing \cite{taft}, which may lead to potential instabilities. 

In this paper, we propose an LEM for energy transactions at the distribution level. The LEM consists of a two-tier structure. The lower level consists of DER-coordinated assets (DCA) located at each secondary feeder bidding into a secondary level market. DCAs are entities that aggregate and coordinate the DERs within their secondary feeder in order to bid into the SM. These DERs could include renewable generation such as rooftop solar PV as well as battery storage and/or flexible loads. We note here that our market mechanism does not inherently rely on any assumptions about ownership structures, e.g. it could be possible for a single agent to coordinate DCAs across multiple primary nodes. An SM operator (SMO) is assumed to oversee the market operations at this level, clearing and scheduling the DCAs. At the upper level, the SMO in turn bids into a PM as an agent representing a primary feeder node. These bids are in turn cleared and scheduled by a PM Operator (PMO), who represents a primary feeder. The payments made by the PMO to the agents at the primary feeder nodes, i.e. the SMOs, are denoted as d-LMPs, and those that are made in turn by the SMO to the DCAs are denoted as local retail tariffs. 

Both the lower and the upper level market solutions proposed here are based on an optimization framework, with the upper level based on a distributed Proximal Atomic Coordination (PAC) approach \cite{rabab_tsg,romvary2021proximal}, while the lower level uses decentralized optimization at each primary feeder node. While the upper level accommodates detailed power physics including nonlinear DistFlow based power balance, and various capacity limits on the main decision variables, the lower level accommodates accurate forecasting of generation and consumption of various DERs with finer granularity and therefore better accuracy. The lower level market also incorporates oversight over the DERs’ actual participation in the market and any unmet commitments thereof in the form of a commitment score. Suitable accommodation is made in recognizing and reflecting any vulnerabilities that may be present in the form of security breaches at the secondary feeder level. Both the PMO and SMO are proposed to be managed by a DSO. We note that there may be multiple primary feeders and thus multiple PMOs connected to a single substation. In this case, all these PMOs would be coordinated by the DSO at the substation. However, for brevity, we assume in this paper that there is only one primary feeder (and PMO) per substation. The overall structure of this hierarchical market is illustrated in \cref{fig:schematic}.

Overall the incorporation of such a hierarchical market framework for DERs allows an efficient incorporation of various expanded responsibilities in a local market. A DSO needs to take on various roles, including maintaining system reliability, facilitating transactions between agents and aggregators, and enabling energy procurement, market clearing and scheduling. Our proposed PMO-SMO structure distributes these roles between the two tiers with greater emphasis on grid physics in the upper level and addressing consumer preferences, reliable performance of DERs and monitoring of security breaches in the lower level. Through a modified case study of an IEEE-123 bus primary feeder test case with multiple secondary feeders at each bus, we demonstrate the functioning of the hierarchical structure and show that the LEM can coordinate and aggregate local DERs more effectively, and enable an optimal combination of local power and power drawn from the bulk grid. This in turn helps reduce distribution level costs and d-LMPs. The incorporation of a commitment score helps to maintain better reliability while still extracting flexibility from customers and DERs. Finally, the time-varying local retail tariffs lead to more efficient market scheduling and lower final costs for end-users, while ensuring that DERs and consumers are correctly compensated for the flexibility services they provide to the grid.

\subsection{Related work}

Several papers have addressed the topic of LEMs and can be grouped into three broad categories – (i) local markets, (ii) hierarchical market structures, and (iii) real-world deployments of local electricity markets. 

Category (i) deals with papers that introduce the concept of LEMs and related solutions \cite{hvelplund2006renewable,kim_bell_labs_microgrid,chen2018next,bray2018unlocking,nemneb}. Reference \cite{hvelplund2006renewable} is the earliest reference in the literature for the term ``local electricity market''. Reference \cite{kim_bell_labs_microgrid} describes how microgrids provided a way of aggregating smaller resources to participate in a market structure. Reference \cite{chen2018next} deals with the concept of having customers and smaller sized DERs participating in a market structure. Transactive energy is also a big driver in enabling LEMs, as compensating consumer resources for services rendered is a key concept in LEMs. LEMs can in fact be viewed as specific structural realizations of a transactive energy framework, and consider the wider system impact and interaction with the Wholesale Energy Market (WEM) \cite{bray2018unlocking}. Other methodologies exist which enable or incentivize the participation of DERs in the WEM, such as net energy metering (NEM), and net energy billing (NEB) \cite{nemneb}. However, these solutions do not encourage full participation of the resources and are often restrictive in their implementation \cite{adapen}. 

Category (ii) corresponds to papers that layout LEMs with a  direct interconnect to the WEM \cite{lezema2019lem,olivella2018optimization,manshadi2015hierarchical,rabab_tsg,chen2018next,bjarghov2021developments,sousa2019peer}. This differs somewhat from papers in Category (i) wherein the local market structures were largely theoretical constructs and still evolving with standardizations yet to emerge. Reference \cite{lezema2019lem} details the interaction between WEM and LEM, and provides numerical results on the cost savings provided by LEM. However, the paper does not disaggregate the price at the DSO and the consumer level, it rather uses a uniform price throughout. References \cite{olivella2018optimization,manshadi2015hierarchical} propose alternative retail market structures that interact seamlessly with a WEM. They utilize a centralized optimization framework, an objective function which aims to minimize the operational costs for the market operator or the DSO, and primarily consider the market participants to be microgrids and/or aggregators. In contrast to these references, our earlier work in \cite{rabab_tsg} proposed a distributed optimization framework which was used to minimize a combination of social welfare and line losses, subject to OPF based on nonlinear Distflow. In this framework, general agents representing DERs can bid into a local market at the primary feeder level, which interacts directly with a WEM. In addition to these papers, several surveys that capture the evolving LEM landscape have been carried out \cite{chen2018next,bjarghov2021developments,sousa2019peer}. Of these, in \cite{chen2018next}, the authors lay out the evolving market structure which will enable customer participation in a market structure. Reference \cite{bjarghov2021developments} is a survey that details the related work in LEM design, existing theoretical tools and models studied in the context of LEM, and challenges of realizing an LEM structure. Reference \cite{sousa2019peer} carries out a survey of peer-to-peer markets. Our paper is similar to those considered in \cite{bjarghov2021developments}, and is a distinct addition in the form of a two-level, hierarchical LEM that includes grid physics, accurate forecasts of DER generation, DER characteristics such as follow-through or unmet commitment, and vulnerability to security breaches.

Category (iii) corresponds to works related to real-field implementation \cite{nyutilitydive,mengelkamp2018designing,thomas2020local,luth2020distributional}. A pioneer solution to integrate DERs has been implemented in the New York electricity market, in the United States \cite{nyutilitydive}. This reference as well as \cite{mengelkamp2018designing} underscore the potential that DER participation brings in terms of coordinated system operation, and shows the feasibility of smaller DER agents effectively participating in electricity markets. In addition to such efforts in the United States, there are several illustrations of successful LEMs in Europe \cite{thomas2020local,luth2020distributional}. Reference \cite{thomas2020local} proposes a two-stage auction based local market mechanism to allocate physical storage rights. Reference \cite{luth2020distributional} discusses the recently proposed market design rules in the current context of the German market with numerical simulations, and a novel market design called Tech4all is introduced. The book \cite{lembook} is an excellent reference for an overall state of the art summary on LEM.

\subsection{Our approach}
Our paper proposes an LEM that connects with the WEM, and belongs to Category (ii). The LEM consists of a two-level hierarchical structure, located in the distribution grid, consisting of a secondary market (SM) at the lower level and a primary market (PM) at the upper level. Together, this hierarchical structure allows for an efficient functioning of the distribution grid which has to achieve multiple objectives and satisfy complicated constraints, by accommodating complex grid physics such as nonlinear and possibly unbalanced power flows, at the upper level by the PM and consumer needs and constraints, in the SM. 

The PM and SM are assumed to be operated by an SMO and a PMO, respectively.  The SM consists of DCA bids submitted to the SMO who clears and schedules them to determine local electricity tariffs. The underlying optimization framework determines the schedules of real and reactive power injections of these DCAs (located at each secondary feeder) as well as their optimal flexibility ranges. Also included are constraints corresponding to bid flexibilities, budget, power balance and capacity limits. A multiobjective cost function of disutility to the DCAs, net cost to the SMO, bid-commitment reliability, and bid-flexibility, is utilized. 

The PM consists of bids from each SMO at node $i$ in the primary feeder, which are aggregations of the cleared market schedules of all the secondary feeders connected to $i$. The objective function includes a weighted combination of social welfare and line losses. The constraints on real and reactive power injections are determined using the net flexibility range of each SMO from the SM clearing. Power balance constraints based on nonlinear DistFlow, capacity limits on $P$ and $Q$ injections,  thermal line limits and bounds on voltages are included. The underlying optimization framework is based on the distributed approach in \cite{rabab_tsg}, and employs a PAC algorithm to determine the schedules of these agents at each $i$. The resulting market clearings include the real and reactive power injections at each primary feeder node and the d-LMPs which correspond to the payments made by the PMO to the SMOs (or vice versa). The net injection from the entire primary feeder is conveyed at the distribution substation to the WEM.

Together with the SM at the lower-level and PM at the upper level, our LEM is used to determine local electricity tariffs and d-LMPs that correspond to the payments between SMOs and DCAs, and between PMOs and SMOs, respectively. Both these prices capture the fine-grain locational and temporal variation in a distribution grid, and form the basis of an efficient LEM that can efficiently allocate resources and accurately compensate prosumers for grid services provided. In this paper, we primarily focus on flexibility services that DERs and flexible loads can offer. This flexibility could be achieved by DCAs through several different types of actions, some examples include (i) load shifting or curtailment, (ii) dispatching distributed generation (such as diesel generators) or battery storage, (iii) using smart inverters to curtail active power from non-dispatchable renewable resources like rooftop PV, or alter $Q$ injections by varying the power factor. In addition to determining the prices, the LEM is also used to determine and update commitment scores at the secondary feeder level, which quantify the ability of each DCA to fulfil their contractual commitments and follow the cleared market schedules. We validate the entire LEM using a modified IEEE 123-node test network, with 14\% renewable generation in terms of nameplate capacity\footnote{relative to a peak load of 3.6 MW on the entire feeder} and \textit{up to} 50\% flexible consumption\footnote{relative to the baseline or nominal load} distributed over 79 primary feeder nodes, and between 3-5 secondary feeders at each of the nodes.  

\subsection{Our contributions}

We propose a novel local electricity market for real-time energy transactions in a distribution grid with high DER penetration. The following are its unique features and the key contributions of our paper:
\begin{itemize}
    \item Hierarchical local electricity market structure that is electrically collocated with the current distribution network.
    \item Effectively address multiple functions of the distribution grid by virtue of the proposed LEM's hierarchy, with grid physics considered in the PM and consumer needs and constraints in the SM.
    \item Systematic approach showing how DERs, along with their flexibility bids, can be coordinated and aggregated in real-time via DCAs, and how these aggregated entities can participate in retail and/or wholesale electricity markets.
    \item Optimization of multiple objectives including commitment reliability, net cost, flexibility and utility in the SM; and net costs, utility and line losses in the PM. 
    \item Accommodation of disparate constraints including budget, power balance, capacity and flexibility limits in the SM; and nonlinear ACOPF constraints in the PM.
    \item Generation of a novel commitment score aimed at improving the reliability of our LEM by tracking the performance of DCAs at the secondary level.
    \item Fine-grain pricing in the form of local retail tariffs in the SM and d-LMPs in the PM that vary with both location and time, allowing for more efficient allocations in terms of lower costs and accurate compensation of DERs. This in turn provides an alternative to the current practices of net-metering and/or direct participation of DERs through aggregators at the WEM level (under FERC 2222).
    \item Validation of the entire LEM using a modified IEEE 123-bus, with a high penetration of DERs and flexible loads.
\end{itemize}

The rest of the paper is organized as follows. In \cref{sec:methodology} we introduce the structure of the LEM, including the lower level SMO, upper level PMO, and the interactions between the PMO, SMO, and WEM. In \cref{sec:results} we present numerical results on a modified IEEE-123 test feeder, with high levels of PV penetration and load flexibility. In \cref{sec:conclusion} we provide concluding remarks. 

\section{A local electricity market (LEM) \label{sec:methodology}}

\begin{figure}
    \centering
    \includegraphics[scale=0.5]{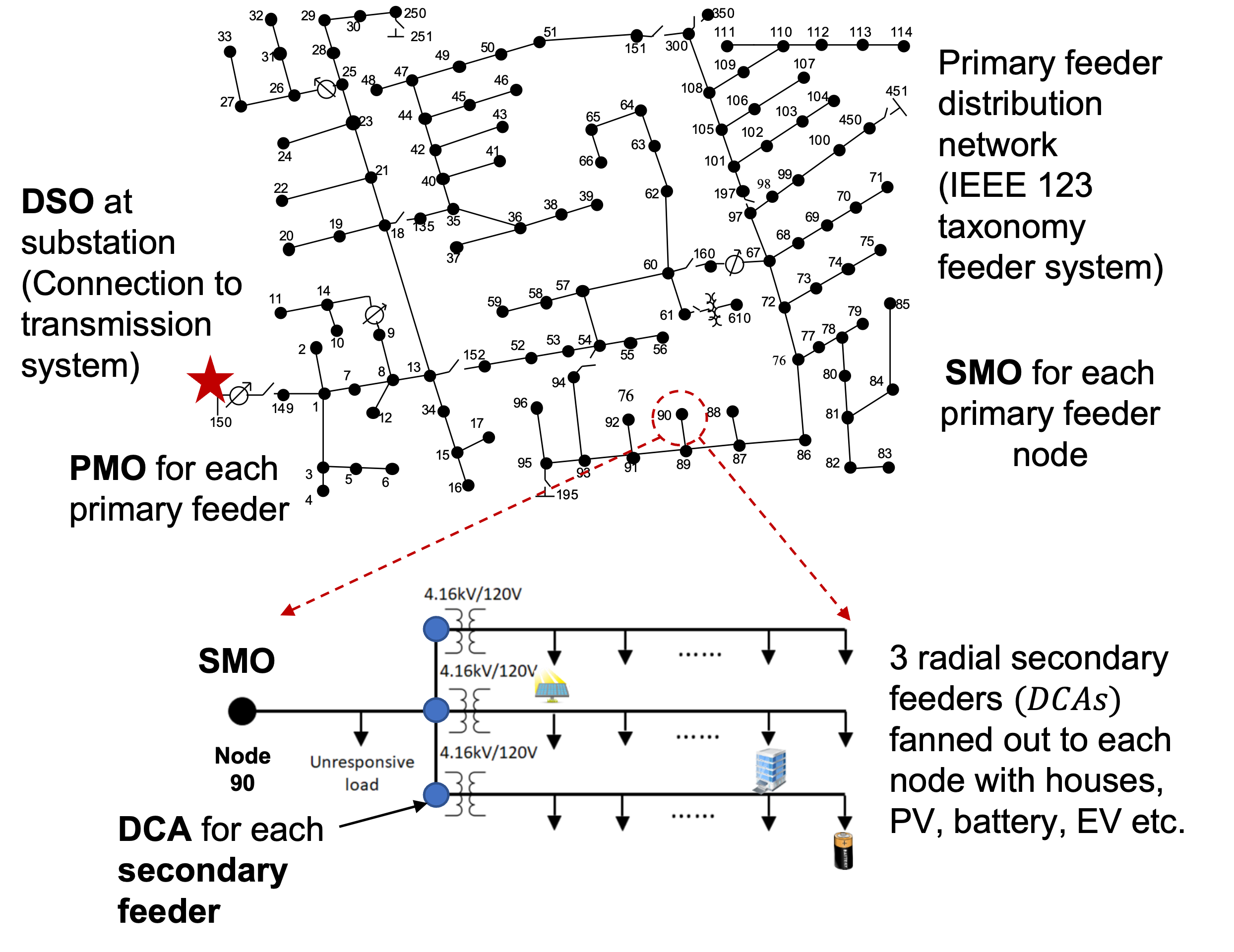}
    \caption{Overall schematic illustrating how the proposed LEM integrates seamlessly into the existing radial distribution network, and connects with bulk transmission.}
    \label{fig:schematic}
\end{figure}
\begin{figure}
    \centering
    \includegraphics[scale=0.5]{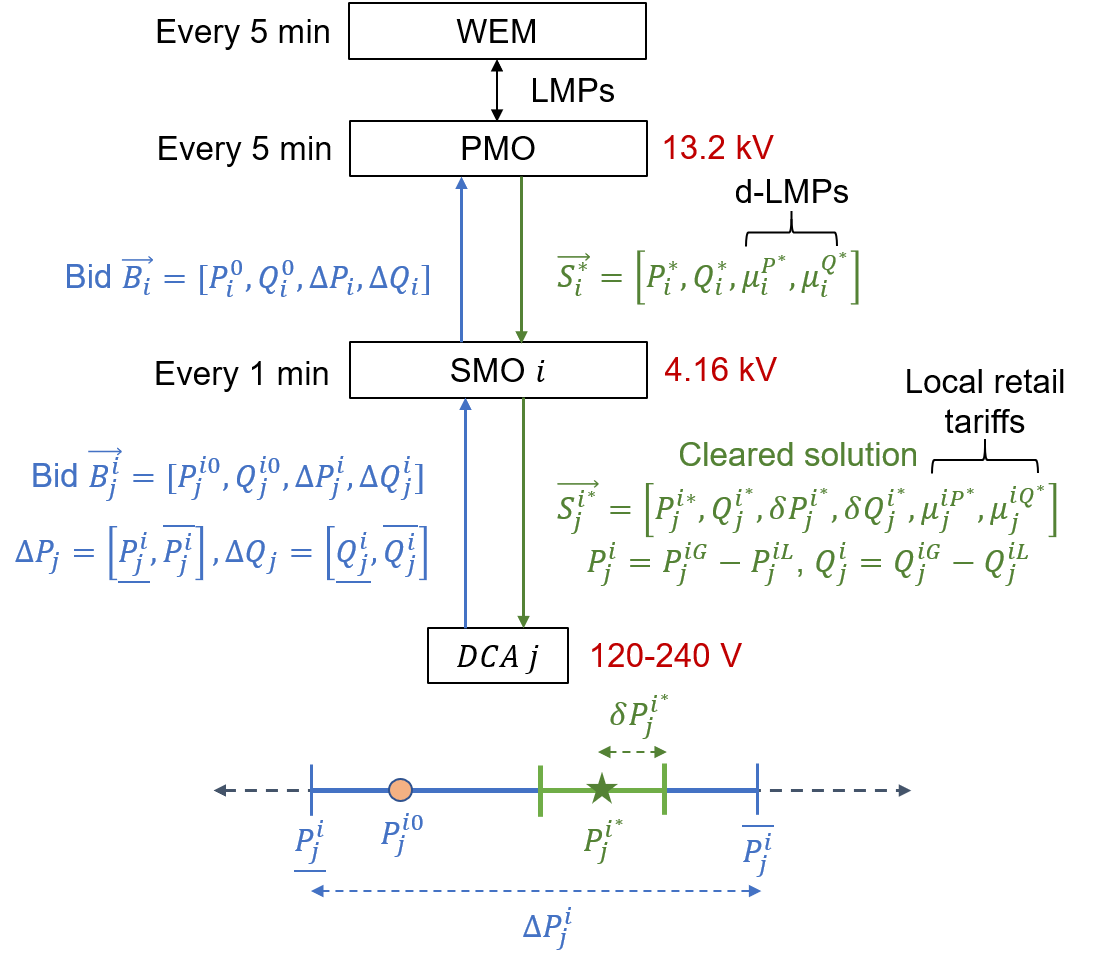}
    \caption{Proposed hierarchical LEM structure and optimization framework incorporating bid flexibility. The upper level PMO coordinates with the WEM, while the SMO at the lower level oversees the DCAs. At the bottom of this figure, we show an example of a feasible flexibility bid from a DCA and the revised flexibility range that results from the SM clearing.}
    \label{fig:structure}
\end{figure}

The LEM we propose in this paper consists of a two-level hierarchical structure, located in the distribution  grid, with an SM at the lower level and a PM at the  upper level, a schematic of which is shown in \cref{fig:schematic,fig:structure}. These markets are assumed to be operated by an SMO and a PMO, respectively, with combined oversight of both provided by a DSO. In what follows, we provide details of the SM and PM, the interface between the two, and their timelines. Throughout this paper, we define net injections as generation less load (i.e. $P = P^G - P^L, Q = Q^G - Q^L$), thus net generation would be positive while net loads are taken to be negative. All bids are assumed to be based on load/generation forecasts one timestep into the future, for both the PM and SM. Thus, the cleared schedules hold for the entire duration of the next period (i.e., either $\Delta t_p$ or $\Delta t_s$), until the next time market bidding and clearing occur. 

\subsection{SM structure}

The following quantities are defined before we pose the optimization problem that serves as the backbone of the SM:
\begin{itemize}
        \item $i \in \mathcal{N}_I$: Set of indices of all SMOs under consideration, downstream of a particular PMO.
        \item $j \in \mathcal{N}_{J,i}$: Set of indices of all DCAs under consideration, downstream of a given SMO $i$. 
        \item $P_j^{i0}, Q_j^{i0}$: Baseline net active and reactive power injections for DCA $j$, which could either be forecasted by the SMO based on historical data or explicitly submitted by the DCA as part of its bid. 
        \item $\Delta P^i_j = [\underline{P}^i_j,\overline{P}^i_j]$, $\Delta Q^i_j = [\underline{Q}^i_j,\overline{Q}^i_j]$: Bid flexibilities for each DCA. These intervals represent the range of maximum downward and upward flexibilities in $P$ and $Q$ injections being offered by the DCA.
        \item $t_p$ and $t_s$: Timestamps for the PM and SM, respectively.
        \item $\Delta t_p, \Delta t_s$ and $\Delta t_{WEM}$: Time periods for the PM, SM, and WEM, respectively.
        \item $\hat t_p$: The most recent PM clearing prior to the current SM interval [$t_s,t_s+\Delta t_s$].
        \item $P_i^*(\hat{t}_p)$ and $Q_i^*(\hat{t}_p)$: Real and reactive power setpoints, respectively, provided by the PM to SMOs at time $\hat{t}_p$.
        \item $n_p, n_s$: Number of primary clearing periods and number of secondary clearings per primary period, respectively.
\end{itemize}
The following decision variables are determined as outputs of the optimization for each DCA $j$ bidding to SMO $i$: 
\begin{itemize}
    \item Decision vector $\vec{S}_j^i = [P_j^i,Q_j^i,\delta P_j^i, \delta Q_j^i,\mu_j^{iP},\mu_j^{iQ}]$
    \item $P_j^i,Q_j^i$: Optimal power injections assigned as setpoints by the SMO $i$ to DCA $j$.
    \item $[\delta P_j^i, \delta Q_j^i]$: Optimal symmetric flexibility ranges around the above setpoints $P_j^i, Q_j^i$, i.e. the DCA is directed to have net injections within these intervals $[P_j^i - \delta P_j^i, P_j^i + \delta P_j^i]$, $[Q_j^i - \delta Q_j^i, Q_j^i + \delta Q_j^i]$.
    \item $\mu_j^{iP},\mu_j^{iQ}$: Local electricity tariffs.
    \item $C_j^i(t)$: Commitment score where $ C_j^i(t) \in [0,1]$ reflects the SMO's confidence in whether the DCA $j$ will reliably follow their committed injections within the flexibility range specified above.
\end{itemize}

\noindent The solutions correspond to power injections and net tariffs corresponding to each DCA, determined at $t_s$ and applied over the period $[t_s, t_s+\Delta t_s]$, for all $t_s$. The following objective function and constraints define the underlying multi-objective optimization problem $\forall \; t_s, \; t_p \; (t_p \leq t_s)$. All quantities and variables in \cref{eq:opt} are specified for the current secondary timestep $t_s$ unless explicitly specified otherwise.
\begin{subequations}
\label{eq:opt}
\begin{align}
& \min_{\vec{S}_j^i} \sum_{j \in \mathcal{N}_{J,i}} \{f_{j,1}^i,f_{j,2}^i,f_{j,3}^i,f_{j,4}^i\} \label{eq:cost} \\
& f_{1,j}^i \succ f_{2,j}^i \succ f_{3,j}^i \succ f_{4,j}^i \label{eq:ranking} \\ 
& f_{j,1} = -C^i_j (P_j^{i^2} +  Q_j^{i^2}), \; f_{j,2}^i  = \mu_j^{iP} P_j^i + \mu_j^{iQ} Q_j^i \nonumber \\
& f_{j,3}^i = -(\delta P_j^i + \delta Q_j^i) \nonumber  \\
& f_{j,4}^i = \beta_j^{iP}(P_j^i - P_j^{i0})^2 + \beta_j^{iQ}(Q_j^i - Q_j^{i0})^2 \nonumber \\
& \text{subject to:} \nonumber\\
& \; \underline{P}_j^i + \delta P_j^i \leq P_j^i \leq \overline{P}_j^i - \delta P_j^i \label{eq:Plim} \\ 
& \underline{Q}_j^i + \delta Q_j^i\leq Q_j^i \leq \overline{Q}_j^i - \delta Q_j^i \label{eq:Qlim} \\ 
& \delta P_j^i, \; \delta Q_j^i \geq 0, \label{eq:flex} \\
& 0 \leq \mu_j^{iP} \leq \overline{\mu}^{iP}, 0 \leq \mu_j^{iP} \leq \overline{\mu}^{iQ}  \; \label{eq:fairness} \\ 
& \sum_{t_s}^{t_s + \Delta t_p} \sum_{j \in \mathcal{N}_{J,i}} \mu_j^{iP}(t) P_j^i(t) \Delta t_s \leq \mu_i^{P^*}(\hat{t}_p) P_i^*(\hat{t}_p) \Delta t_p \label{eq:budgetP} \\
& \sum_{t_s}^{t_s + \Delta t_p} \sum_{j \in \mathcal{N}_{J,i}} \mu_j^{iQ}(t) Q_j^i(t) \Delta t_s \leq \mu_i^{Q^*}(\hat{t}_p) Q_i^*(\hat{t}_p) \Delta t_p \label{eq:budgetQ}\\ 
& \sum_{j \in \mathcal{N}_{J,i}} P_j^i(t_s) = P_i^*(\hat{t}_p), \quad \sum_{j \in \mathcal{N}_{J,i}} Q_j^i(t_s) = Q_i^*(\hat{t}_p) \label{eq:PQbalance}
\end{align}
\end{subequations}
The cost functions in~\eqref{eq:cost} correspond to the following: 
\begin{enumerate}
    \item \textbf{Commitment $\mathbf{f^i_{j,1}}$}: This term maximizes the injections assigned to more trustworthy DCAs (i.e., $C^i_j$ closer to 1) while minimizing the scheduling of DERs with lower commitment scores, who are relatively less likely to follow through on their contractual commitments.
    \item \textbf{Net costs $\mathbf{f^i_{j,2}}$}: This term minimizes the \textit{net} costs to the SMO for running its SM, which are composed of payments made out by the SMO to DCAs that are net generators, and denote revenue if the DCAs are loads. 
    \item \textbf{Flexibility $\mathbf{f^i_{j,3}}$}: These terms aim to maximize the aggregate flexibility that the SMO can extract from its DCAs, and in turn offer to the PMO.
    \item \textbf{Disutility $\mathbf{f^i_{j,4}}$}: These terms aim to minimize the disutility or inconvenience caused to DCAs when they provide flexibility to the operator. Thus, our SMO is an altruistic entity that also considers welfare maximization for its DCAs. For our simulations in \cref{sec:sec_mkt}, the disutility coefficients were chosen as $\beta^{iP}_j, \beta^{iQ}_j \sim \mathcal{U}[0.1,1]$.
\end{enumerate}

The multiobjective optimization problem in \cref{eq:opt} was formulated and solved using a hierarchical or lexicographic optimization-based approach. This method has been widely used in the literature to solve multiobjective problems \cite{marler_survey_2004, gunantara2018review}, particularly when the objectives have different units and may not be comparable in magnitude. In the hierarchical method, the different objectives are ranked in descending order in terms of their importance to the decision maker. The SMO orders their four objectives as shown in \cref{eq:ranking}, assigning commitment reliability as the most important goal and DCA disutility as being the least important. The SMO then solves a series of optimization problems, sequentially optimizing each of these objectives one at a time, in descending order of importance:
\begin{align}
    \min_{\vec{S}_j^i} & \; F_k = \sum_{j \in \mathcal{N}_{J,i}} f_{j,k}^i (\vec{S}_j^{i}) \; \forall \; k = 1, 2, 3, 4 \\
    \text{s.t.} \; & f_{j,\ell}^i(\vec{S}_j^{i}) \leq (1 + \epsilon)\sum_{j \in \mathcal{N}_{J,i}}  f_{j,\ell}^i (\vec{S}_j^{i^*}) = (1 + \epsilon)F_{\ell}^*, \label{eq:degrade} \\
    & \forall \; \ell = 1,2,\dots, k-1, \; k>1 \nonumber \\
    & \text{constraints in eqs.} \; (\ref{eq:Plim}) \; \text{to} \; (\ref{eq:PQbalance})
\end{align}
At each step after optimizing the 1st objective (commitment reliability $f^i_{j,1}$), additional constraints are placed on the values of the previously optimized objectives as in \cref{eq:degrade}. The hyperparameter $\epsilon$ controls the extent to which previous objective values are allowed to be degraded while searching for the new minima. We used $\epsilon = 0.05 \; (5\%)$ for our simulations. Using such a hierarchical approach also helps us get around the issues of our objective terms potentially being on different orders of magnitudes since we only minimize a single objective at each step, and the constraints in \cref{eq:degrade} look at the \textit{relative} changes in the objective function values rather than their absolute magnitudes. Thus, we can proceed without needing to normalize any of the terms. This is especially advantageous here since most normalization methods require either prior knowledge of maximum or minimum objective values (which are not known beforehand in our case), or entail solving additional optimization problems to find these values at every iteration, which can be computationally expensive \cite{GrodzevichNormalizationOptimization,Makler2006Function-transformationOptimization, Miettinen2002OnOptimization}. Finally, we note that multiobjective optimizations problems in general do not have unique minima. Rather, the goal here is to find a \textit{Pareto}-optimal set or efficient frontier of multiple possible solutions \cite{marler_survey_2004,Deb2014Multi-objectiveOptimization}.

The constraints in \cref{eq:Plim,eq:Qlim,eq:flex} reflect the feasible flexibilities for each of the DCAs. Constraints in \cref{eq:fairness} enforce that the cleared tariffs $\mu_j^i$ for any DCA $j$ cannot exceed a price ceiling $\overline{\mu}^i_j$. These upper bounds may be set by an external regulatory authority such as a DSO, Independent System Operator (ISO), Regional Transmission Operator (RTO), or Public Utility Commission (PUC). The budget constraints in \cref{eq:budgetP} and \cref{eq:budgetQ} ensure that the total net payments made out by the SMO to its DCAs over all the SM clearings within each primary interval $[t_s,t_s+\Delta t_p]$ do not exceed its net revenue received from the PMO during the same period. We enforce these budget constraints separately for $P$ and $Q$ since proposed reactive power markets often behave quite differently compared to the energy market for real power \cite{safari_reactive_2020}. Note that our optimization problem solves for the prices in terms of [\$/kW] or [\$/kVAR] so these are converted to [\$/kWh] or [\$/kVARh] respectively, before applying any of the budget constraints.

The constraints as written in \cref{eq:budgetP} and \cref{eq:budgetQ} can be quite restrictive since they require the SMO to balance its budget for \textit{every} PM period $\Delta t_p$. Thus, in addition to the strict budget constraints above, we also considered a softer implementation of the same. These require the SMO to balance its budget over the course of a longer time horizon (e.g. 1 day), allowing the SMO to run a deficit for some PM clearings, if needed. These are shown in \cref{eq:budgetPrelax,eq:budgetQrelax}, where we consider a planning horizon of $n_p$ primary clearing periods for the budget balance:

\begin{align}
    & \sum_{0}^{n_p\Delta t_p} \sum_{j \in \mathcal{N}_{J,i}} \mu_j^{iP}(t) P_j^i(t) \Delta t_s \leq \sum_{0}^{n_p\Delta t_p} \mu_i^{P^*}(t) P_i^*(t) \Delta t_p \label{eq:budgetPrelax} \\
    & \sum_{0}^{n_p\Delta t_p} \sum_{j \in \mathcal{N}_{J,i}} \mu_j^{iQ}(t) Q_j^i(t) \Delta t_s \leq \sum_{0}^{n_p\Delta t_p} \mu_i^{Q^*}(t) Q_i^*(t)  \Delta t_p \label{eq:budgetQrelax}
\end{align}

\noindent Since the optimization problem in \cref{eq:opt} is currently framed as a single period problem, we transformed our inherently multiperiod budget constraints to a \textit{quasi-multiperiod} form for implementation. This was done by assuming that the SMO evenly redistributes its leftover net revenue over all the remaining SM clearings in the current budget period. For example, the quasi-multiperiod version of the strict budget constraint in \cref{eq:budgetP} is given by: 
\begin{align}
    & \sum_{j \in \mathcal{N}_{J,i}} \mu_j^{iP}(t_s) P_j^i(t_s) \Delta t_s \leq \label{eq:quasimp} \\ 
    & \frac{\sum_{0}^{\hat{t}_p} \mu^{P^*}_i(t) P^*_i(t)\Delta t_p - \sum_{0}^{t_s - \Delta t_s} \sum_{j \in \mathcal{N}_{J,i}} \mu^{iP^*}_j(t) P^{i^*}_j(t) \Delta t_s}{\frac{\hat{t}_p + \Delta t_p - t_s}{\Delta t_s}} \nonumber
\end{align}
where $\hat{t}_p$ is the most recent primary clearing time  prior to this SM clearing, and the denominator is the number of secondary clearings left in the current PM interval $[\hat{t}_p, \hat{t}_p + \Delta t_p]$. A similar quasi-multiperiod form can be derived for the relaxed budget constraints in \cref{eq:budgetPrelax,eq:budgetQrelax} as well.

Finally, \cref{eq:PQbalance} denotes power balance constraints for the SMO, where the sum total of net injections from the DCAs downstream needs to satisfy the net flows from the primary feeder upstream, These net injections $P_i^*(\hat{t}_p)$ and $Q_i^*(\hat{t}_p)$, at each primary feeder node (SMO) $i$ are scheduled by the PMO. Since the PM clears less often than the SM, these values can be treated as constant for the SM optimization problem over each $\Delta t_p$. Suitable convex relaxations were added to ensure feasibility and minimal optimality gap, and the resulting optimization problems were solved using Gurobi in Python\footnote{\url{https://www.gurobi.com/}}. 


The overall operation of the SM is summarized as follows: Starting with DCA bids $\vec{B}_j^i = [P_j^{i0}, Q_j^{i0}, \Delta P^i_j, \Delta Q^i_j]$ as inputs, the constrained optimization problem in eqs. (\ref{eq:cost})-(\ref{eq:PQbalance}) is solved by the SMO, with the market clearing resulting in net-injections $P_j^{i^*}(t_s)$ and local electricity tariff $\mu^{iP^*}_j(t_s)$, for real power. A similar set of injections and tariffs are derived for reactive power as well. Each of these solutions corresponds to the net-injection at DCA $j$, cleared at $t_s$, and holds for the following period $[t_s,t_s+\Delta t_s]$. Together with these net injections, the SMO also obtains, as a part of the above optimization procedure, a feasible flexibility $\delta P_j^i$ in the real power injection for node $j$, and a similar flexibility for the reactive power injection. Together, the complete set of solutions from the SMO at secondary feeder $j$ connected to the PMO at node $i$ is given by $S_j^{i^*}$ (see \cref{fig:structure}). These market cleared solutions are used to establish bilateral contracts between the SMO and its DCAs, with localized retail tariffs for the power injections, differentiated for each of them. 


Another unique aspect of this SM structure is the commitment score $C_j^i(t)$. This is determined at every $t_s$, with a score decrease following every event where the DCA exhibits unmet commitment. As the goal is to have the SMO $i$ reduce this score for such events and to reward the DCA $j$ when they do follow through on their commitment, we propose the following recursive algorithm for a continuous update of $C_j^i(t_s)$:

\begin{align} 
    C^i_j(t_s) & = C^i_j(t_s - 1) - \frac{\widetilde{e^{iP}_j} (t_s) + \widetilde{e^{iQ}_j} (t_s)}{2} \; \forall \; j \in \mathcal{N}_{J,i} \label{eq:commit3} \\
    e^{P}_j(t_s) & = \llbracket\hat{P}_j > \overline{P}_j^{i^*}\rrbracket (\hat{P}_j - \overline{P}_j^{i^*}) + \llbracket\hat{P}_j < \underline{P}_j^{i^*}\rrbracket(\underline{P}_j^{i*} - \hat{P}_j) \nonumber\\
    & + \llbracket\underline{P}_j^{i^*} \leq \hat{P}_j \leq \overline{P}_j^{i^*}\rrbracket\max (\hat{P}_j - \overline{P}_j^{i^*},\underline{P}_j^{i^*} - \hat{P}_j) \label{eq:commit1}\\
    \widetilde{e^{iP}_j} (t_s) & = \frac{e^{iP}_j(t_s)}{|P_j^{i^*}(t_s)|}, \quad \widetilde{e^{iQ}_j} (t_s) = \frac{e^{iQ}_j(t_s)}{|Q_j^{i^*}(t_s)|} \label{eq:commit2}\\
    \widetilde{e_i^P}(t_s) & = \frac{\mathbf{e_i^P(t_s)}}{{\|\mathbf{e_i^P}}(t_s)\|}, \; \widetilde{e_i^Q}(t_s) = \frac{\mathbf{e_i^Q}(t_s)}{\|\mathbf{e_i^Q}(t_s)\|} \label{eq:commit4}
\end{align}
where $\underline{P}_j^{i^*} = P_j^{i^*} - \delta P_j^{i^*}, \; \overline{P}_j^{i^*} = P_j^{i^*} + \delta P_j^{i^*}$. Here, the deviations in DCA injections are first normalized by their true setpoints in \cref{eq:commit2}, followed by $L_2$ normalization in \cref{eq:commit4}. These normalizations allow us to assess the relative performance of all DCAs under SMO $i$. The resulting commitment score serves as a metric of DCA reliability, and is proposed as a stand-alone indicator of DER commitment or a lack thereof. Further refinements of the SM could include the  use of $C_j^i$ for determining suitable penalties to these DCAs as $C_j^i$ drops below unity and tightly interconnect with a corresponding ancillary market structure. Yet another component that could be coupled with $C_j^i$ is the DCA's vulnerability to cybersecurity breaches. Currently these discussions are omitted, and will be pursued as part of future work. 
\begin{figure}
    \centering
    \includegraphics[scale=0.5]{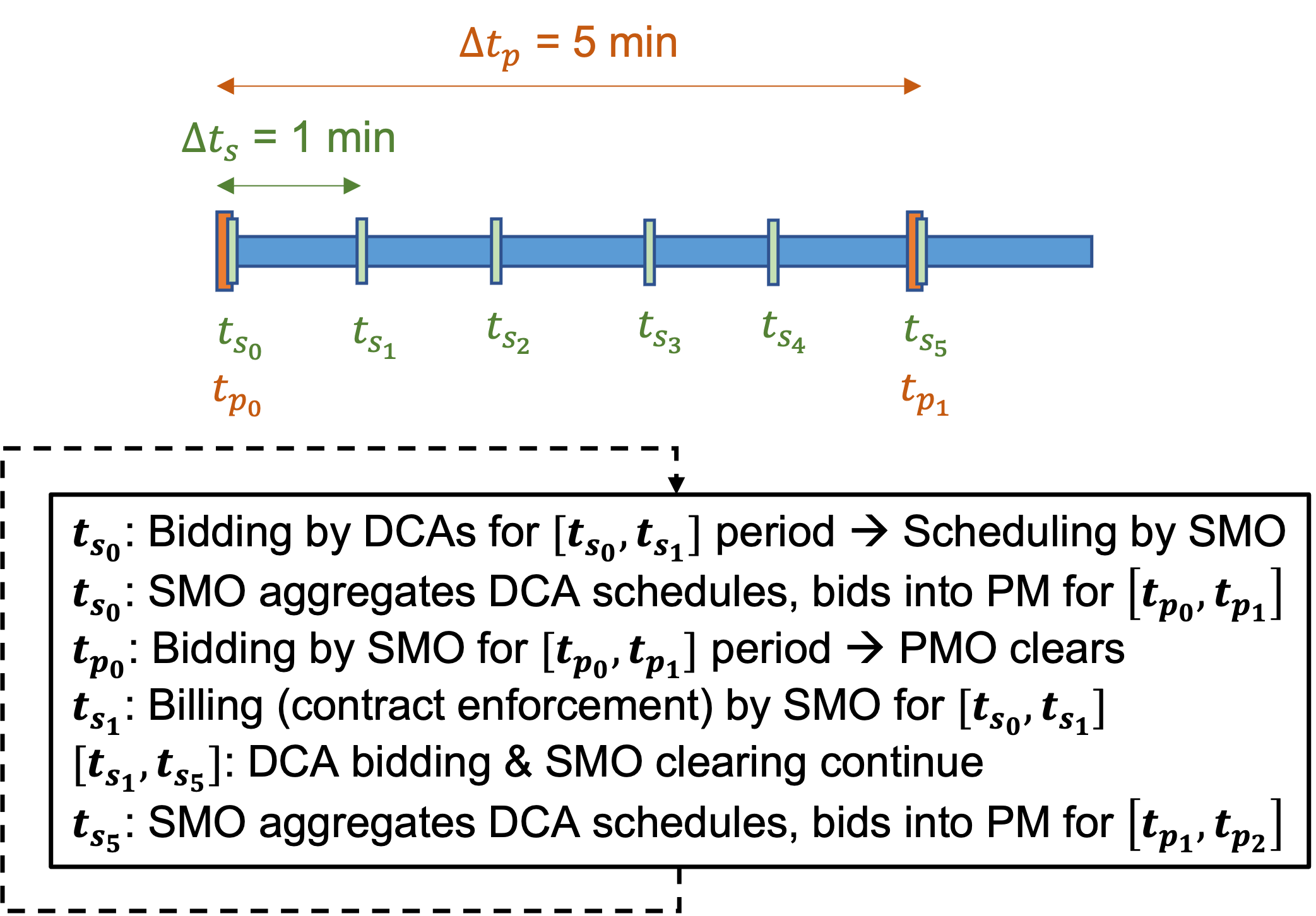}
    \caption{An illustrative timeline of the proposed LEM, also showing the interactions and interplay between the primary and secondary level markets.}
    \label{fig:timeline}
\end{figure}

\subsection{Interface between SM and PM}
\begin{align}
    & P^0_i(t_p) = \sum_{j \in \mathcal{N}_{J,i}}  P_j^{i^*}(t_p),\; Q^0_i(t_p) = \sum_{j \in \mathcal{N}_{J,i}}  Q_j^{i^*}(t_p) \label{eq:aggregation}  \\
    & \Delta P_i = \big[\underline{P}_i = \sum_{j\in \mathcal{N}_{J,i}} P_j^{i^*} - \delta P_j^{i^*}, \overline{P}_i =\sum_{j\in \mathcal{N}_{J,i}} P_j^{i^*} + \delta P_j^{i^*} \big] \nonumber \\
    &  \Delta Q_i = \big[\underline{Q}_i = \sum_{j\in \mathcal{N}_{J,i}} Q_j^{i^*} - \delta Q_j^{i^*}, \overline{Q}_i = \sum_{j\in \mathcal{N}_{J,i}} Q_j^{i^*} + \delta Q_j^{i^*}\big] \nonumber
\end{align}

As mentioned earlier, an SM is located at each primary feeder node $i$, supervised by SMO $i$, and determines market clearings for all secondary feeders $j$ connected to this node with each $j$ assumed to be represented by DCA $j$ in the SM. The market clearing consists of $\left(P_j^{i^*}(t_s), \mu_j^{iP^*}(t_s)\right)$ at each $t_s$, which corresponds to power injections and local electricity tariffs for secondary feeder j. The market also returns an optimal flexibility $\delta P_j^{i^*}$ for each $j$. A similar set of clearing variables are associated with reactive power too. At each primary market timestep $t_p$, this SM clearing is completed \textit{before} the PM clearing at the same timestep. The SMO $i$ in turn uses these DCA solutions from its SM in order to form its own bid into the PM at time $t_p$, representing primary feeder node $i$, that is at the upper level of our proposed LEM. The SMO aggregates power injections from all secondary feeders as its baseline injection $P_i^0$, at time $t_p$. The SMO bid also includes a corresponding aggregated flexibility $\Delta P_i$ based on the DCA-flexibilities $\delta P_j^{i^*}$. The specific bid determinations are given by \cref{eq:aggregation}. These bids and the corresponding PM clearing are assumed to occur every $\Delta t_p$, with the assumption that $\Delta t_p > \Delta t_s$.


\subsection{PM structure}
The starting point for the PM are the SMO bids $P_i$ with the flexibility $\delta P_i$, and a similar set of quantities for reactive power, at each primary feeder node $i$. The PMO (possibly at a substation) which has oversight over the entire primary feeder, accepts these bids, and clears the PM at every $t_p$ at intervals of $\Delta t_p$. The PMO clears the market, with an underlying distributed optimization framework that facilitates market clearing. This framework includes appropriate cost functions such as Social Welfare and line losses, and constraints that correspond to power physics constraints modeled with nonlinear DistFlow equations (branch flow model), as well as various network level constraints. 

The underlying optimization problem for the PMO over a radial distribution grid is defined in \cref{eq:pmo}, where $P_i,Q_i,v_i$ are the nodal real power, reactive power, and voltages, superscripts $G$ and $L$ denote generation and load, and $P_{{i}^\prime,i}, Q_{{i}^\prime,i}, l_{{i}^\prime,i},R_{{i}^\prime,i}, X_{{i}^\prime,i}$ denote the real power, reactive power, squared current magnitude, resistance, and reactance of the branch from node $i^\prime$ to $i$. These nonlinear constraints describe the power flow in the radial distribution grid, assuming balanced flows, small angles, and a convexification of the definition of power using second-order cone constraint programming (SOCP). The remaining constraints describe the network voltage constraints, and generator and load flexibility limits.
\begin{align}
    & \underset{y}{\min} \;\; f^{S-W}(y) \label{eq:pmo} \\
    & \text{subject to:} \nonumber\\
	& v_{i} - v_{{i}^\prime} = \left(R_{{i}^\prime,i}^{2} + X_{{i}^\prime,i}^{2}\right) l_{{i}^\prime,i} - 2 \left(R_{{i}^\prime,i} P_{{i}^\prime,i} + X_{{i}^\prime,i} Q_{{i}^\prime,i}\right)\nonumber\\
	& P_{i}^{G} - P_{i}^{L} = -P_{{i}^\prime,i} + R_{{i}^\prime,i} l_{{i}^\prime,i} + \sum_{k:\left(i,k\right) \in \mathcal{T}} P_{i,k}\nonumber\\
	& Q_{i}^{G} - Q_{i}^{L} = -Q_{{i}^\prime,i} + X_{{i}^\prime,i} l_{{i}^\prime,i} + \sum_{k:\left(i,k\right) \in \mathcal{T}} Q_{i,k}\nonumber\\ 
  	& P_{{i}^\prime,i}^{2} + Q_{{i}^\prime,i}^{2} \leq \overline{S}_{{i}^\prime,i}^{2}\nonumber\\
  	& P_{{i}^\prime,i}^{2} + Q_{{i}^\prime,i}^{2} \leq v_{i} l_{{i}^\prime,i}\nonumber\\
  	& \underline{v}_i \leq v_i \leq \overline{v}_i  \nonumber\\
	& \underline{P}_i^{G} \leq P_i^{G} \leq \overline{P}_i^{G} \nonumber\\
 	& \underline{P}_i^{L} \leq P_i^{L} \leq \overline{P}_i^{L} \nonumber\\
	& \underline{Q}_i^{G} \leq Q_i^{G} \leq \overline{Q}_i^{G} \nonumber\\
  	& \underline{Q}_i^{L} \leq Q_i^{L} \leq \overline{Q}_i^{L} \nonumber
\end{align}

The specific cost function that will be used is as below: 
\begin{align}
    f^{S-W}(y)&=\sum_{j}\Big[f_{i}^{\text{Load-Disutil}}(y) + f_{i}^{\text{Gen-Cost}}(y)\Big]  \nonumber \\
    & + \xi\Big[\sum_{k}\!f_{k}^{\text{Loss}}(y)\Big] \label{eq:cost:1} \\
    f_{i}^{\text{Load-Disutil}}(y)&= \beta_{i}^{P}(P_{i}^{L}-P_{i}^{L0})^2+\beta_{i}^{Q}(Q_{i}-Q_{i}^{L0})^2 \label{eq:cost:2} \\
	f_{i}^{\text{Gen-Cost}}(y)&= \begin{cases}\alpha_{i}^{P} (P_{i}^{G})^2 +\alpha_{i}^{Q} (Q_{i}^{G})^2, & \\
	\lambda_{i}^{P} P_{i}^{G}+\lambda_{i}^{Q} Q_{i}^{G},\text{if } i & \text{ is PCC}
	\end{cases} \label{eq:cost:3}\\
	f_{i^{\prime},i}^{\text{Loss}}(y)&=R_{i^{\prime},i}l_{i^{\prime},i} \label{eq:cost:4}
\end{align}
where $i$ are the indices for SMOs participating in the PM, $P_i^L, Q_i^L$ and $P_i^G, Q_i^G$ are the loads and generation at node $i$ respectively, $P_i^{L0}$ and $Q_i^{L0}$ are the baseline loads, $l_{i^{\prime},i}$ and $R_{i^{\prime},i}$ are the squared magnitude of current and resistance of line from node $i^{\prime}$ to $i$. The cost coefficients for load disutility of the SMO are computed as $\beta_i^{P} = \frac{1}{|\mathcal{N}_{J,i}|} \sum_{j\in \mathcal{N}_{J,i}} \beta^{iP}_j$, where notation $|\mathcal{N}_{J,i}|$ denotes the cardinality of the set $\mathcal{N}_{J,i}$, i.e. the number of DCAs downstream of SMO $i$, and likewise for reactive power disutility $\beta_i^{Q}$. The quadratic coefficients for generating cost are $\alpha_i^P$ and $\alpha_i^Q$, and the wholesale price of power (LMP) from the WEM at the PCC is $\lambda^P$ and $\lambda^Q$. The coefficient of power loss
$\xi$ in \cref{eq:cost:1} reflects the trade off between minimizing line losses versus minimizing disutility and generation costs. These coefficients would be privately chosen by the SMO according to their preferences. For our simulations, we used a constant value of $\xi = 100$ for all SMOs. We chose this value in order to balance the tradeoffs of socioeconomic costs versus line losses. Using a $\xi$ value that’s too low would devalue the line losses term entirely due to relative scaling issues, and a value that's too high is also unrealistic because economic and utilitarian decisions drive the market. During actual implementation, this power loss factor would likely have to be tuned by each SMO over time, based on their operational objectives and above mentioned tradeoffs. Finally, the PMO interfaces with a market operator at the substation and bids into the WEM every $\Delta t_{WEM} = 5 \; min$.

The linear cost term in \cref{eq:cost} for the SM accounts for the retail costs to customers and payments to DGs. However, we use a quadratic cost in \cref{eq:cost:3} for the PM in order to account for the additional costs to the SMO for operating the SM and maintaining its distribution network. These include procuring reserves, adequate storage, and standby or auxiliary generation capacity\footnote{Resources like large-scale battery storage, diesel gensets, natural gas peaker plants, etc.} for contingencies, e.g. to meet shortfalls in case of unmet commitments from its DCAs. These \textit{fixed} and amortized cost coefficients for the SMOs were chosen to be between $\alpha^{P, fixed}_i \in [4,8] \; \$/p.u.$ \cite{lazard2018,camera2019renewable}, with $S_{base} = 1 \; MVA$. In addition to these O\&M costs, the SMO also adjusts these cost-coefficients at every $t_p$ by a small variable amount  $\alpha_i^{P, var}$, based on the weighted average of the retail tariffs across all its DCAs from the previous SM clearing. This variable cost component reflects the average cost of electricity in the SM and thus can also be interpreted as a price bid or offer by the SMO to the PMO, resulting from the SM optimization.  
\begin{align}
    & \alpha^{P, var}_i (t_p) = \frac{\sum_{t_p - n_s \Delta t_s}^{t_p} \sum_{j \in \mathcal{N}_{J,i}} \mu^{jP^*}_i \left| P^{j^*}_i \right|}{\sum_{t_p - n_s \Delta t_s}^{t_p} \sum_{j \in \mathcal{N}_{J,i}}\left| P^{j^*}_i \right|} \\
    &  \alpha^{P}_i (t_p) = \alpha^{P, fixed}_i + \alpha^{P, var}_i (t_p)
\end{align}

With the costs and constraints as above, the OPF problem is solved using CVX in MATLAB\footnote{\url{http://cvxr.com/cvx/}} to carry out the market clearing, which consists of the power injection at each SMO $i$ and d-LMP, the electricity price corresponding to primary node $i$ and determined using the dual variable obtained from the OPF. These clearings are implemented through a distributed optimization algorithm denoted as PAC \cite{rabab_tsg,romvary2021proximal} which consists of peer-to-peer communication between neighboring SMOs in an autonomous manner. This makes the computation more tractable and also reduces communication latencies. Upon reaching an agreement with its neighbours, each SMO enters into a bilateral agreement with its PMO, thereby committing to deliver or consume the decided amount of power, at the d-LMP rate. Any net loads consumed by the SMO will be charged the d-LMP, and equivalently, net generation by an SMO will be remunerated at the d-LMP. All payments will be made to/from the PMO. 

\subsection{Timelines of SM and PM}

The overall time scales of the secondary and primary levels of our proposed retail market in relation to the real-time market in a WEM are indicated in \cref{fig:timeline}. Market clearings of the SM and PM are assumed to be every $\Delta t_s$ and $\Delta t_p$ apart, with $\Delta t_s < \Delta t_p$. For the use-case study, we assume that $\Delta t_p = \Delta t_{WEM}$, i.e., the PM and WEM are cleared together in lockstep. We also assume that $\Delta t_s = 1 \; min$ and $\Delta t_p = 5 \; min$, and that the SM clearing occurs arbitrarily quickly as the complexity and dimensionality of the underlying optimization problem is low. 



\subsection{Assumptions, Observations, and Extensions}
The LEM proposed here has been constructed using a hierarchical structure precisely to address the distinct challenges that a distribution grid poses in comparison to a  transmission grid. This hierarchical structure allows an efficient incorporation of multiple objectives and constraints simultaneously present. Since the SM is closer to the end-user in both location and time, we constructed the SM to be more consumer-centric, with costs and constraints pertaining to consumer flexibility and needs. Since the PM, relatively speaking, has a complex set of physical network topologies, we pay greater attention to the physical costs and constraints in its problem formulation. This allows the DSO, overall, to address the varied roles of reliability challenged by grid physics, and flexibility and granularity in location and time challenged by the presence of disparate consumers with varied needs and constraints. In this study, we assume that the SM is cleared more frequently compared to the PM, in order to quickly accommodate any variations that may occur locally at the DCA level. However, this condition is not necessary to operate our hierarchical LEM, allowing both markets to be in lockstep if need be. 

In this paper, we considered a single period optimization problem solved by the SM at each timestep. However, our market structure can be applied to the multiperiod optimization setting as well. We are currently working on extensions of our model where the SMO and PMO perform multiperiod optimization over a planning horizon into the future, using an approach inspired by optimal control or model predictive control (MPC), similar to \cite{Nudell2022DistributedApproach}. For instance, the core structure of the SM optimization problem would remain the same as in \cref{eq:opt} but would now optimize over several timesteps $T$ into the future, subject to similar constraints as in \cref{eq:Plim}-\cref{eq:PQbalance}:
\begin{align}
    \min_{\vec{S}_j^i(t)} \sum_{t = \tau}^{\tau+T} \sum_{j \in \mathcal{N}_{J,i}} \{f_{1,j}^i(t|\tau), f_{2,j}^i(t|\tau), f_{3,j}^i(t|\tau), f_{4,j}^i(t|\tau)\}
\end{align}
where $(t|\tau)$ denotes predictions or estimates of quantities for future periods $t$ made at time $\tau$. Thus, the decision vector $\{\vec{S}_j^i(t)\}_{t=\tau}^T$ is now higher dimensional since it spans multiple SM periods. In addition to future flexibility bids from DCAs, i.e., ($P_j^{i0}(t|\tau), Q_j^{i0}(\tau), \Delta P^i_j(t|\tau), \Delta Q^i_j(t|\tau)$), the SMO also needs to predict future PM solutions ($\mu^{P^*}_i(\hat{t}_p|\tau), P_i^*(\hat{t}_p|\tau),\mu^{Q^*}_i(\hat{t}_p|\tau), Q_i^*(\hat{t}_p|\tau)$ in order to solve the SM multiperiod optimization problem. Herein lies the main challenge of extending to the multiperiod setting. We are exploring several time series forecasting tools like ARIMA\footnote{Auto Regressive Integrated Moving Average}, exponential smoothing, etc. for this purpose. This optimization would result in both binding spot values that apply for the very next timestep, as well as future values for further into the planning horizon that are non-binding. One of the main benefits of the multiperiod approach is that it allows us to implement more realistic formulations of the budget constraint, without having to resort to assumptions as described in \cref{eq:quasimp}. For example, the multiperiod version of the budget balance for active power would be: 
\begin{align}
    & \sum_{t \in \mathcal{T}_s} \sum_{j \in \mathcal{N}_{J,i}} \mu_j^{iP}(t) P_j^i(t) \Delta t_s
    \leq \sum_{t' \in \mathcal{T}_p} \mu_i^{P^*}(t'|\tau) P_i^*(t'|\tau) \Delta t_p \label{eq:budgetP_mp}
\end{align}
where $\mathcal{T}_s$ and $\mathcal{T}_p$ denote the set of SM and PM clearing times within the current planning horizon $T$, respectively. A multiperiod approach would also allow us to include inter-temporal constraints, in order to better optimize the scheduling of energy storage devices like batteries and EVs as well as thermostatically controlled loads (TCLs) such as heating, ventilation, and air conditioning (HVAC) systems and water heaters \cite{Hao2014AncillarySystems}.

Our focus in this paper has been on real-time energy markets. We have not addressed issues such as settlements and billing, as the relevant discussions will have to necessarily include ancillary markets and reserves to  deal with any unmet commitments and supply-demand imbalances, in real-time. Extensions similar to co-optimization \cite{garcia_dynamic_2017} of these different markets are expected to be possible. Currently, bids into the PM are synthesized from the SM clearings through a simple aggregation process. Advanced game-theoretic approaches such as \cite{kardakos_optimal_2014,somma_optimal_2019} have the potential to generate more intelligent bids, and is a topic for future research. The inclusion of the commitment score in our optimization problem is an effort to address consumer-centric constraints that could lead to unmet commitments. This may be due to a variety of factors such as (i) willful reneging on contracts, (ii) malicious behavior due to system compromise or security breaches, and (iii) changes in environmental or weather conditions. The commitment score could also be potentially used for determining penalties and tariffs, thereby leading to a more efficient market design. Details of this effort are part of our future work as well.

\section{Results and Discussion \label{sec:results}}

\subsection{The Use-case}

The hierarchical LEM proposed in \cref{sec:methodology} is evaluated using a modified IEEE-123 test feeder. A GridLAB-D model\footnote{\url{https://www.gridlabd.org/}} was utilized to simulate this test feeder over the course of a 24 hour period. Rooftop PV (with smart inverters) was assumed to be present at nodes 5, 20, 50, 63, and 94, with a total PV generation capacity of 510.3 kW. This corresponded to a DER (PV) penetration of about 14\%, assuming that the peak load is at about 3.6 MW \cite{kersting1991radial,postigo2017review}. An SMO was assumed to be present at 79 of the primary feeder nodes (i.e. $\left|\mathcal{N}_I\right| = 79$), and that flexible loads were present at all of these nodes with each DCA capable of \textit{up to} $\pm$50\% deviations around their baseline injections. We assumed this maximum flexibility based on past studies forecasting demand response potentials in the US \cite{olsen_demand_2014}. The GridLAB-D model included triplex meters to record P and Q injections every minute, at each of these 79 nodes. Weather data for Boston, MA was used to forecast PV generation, and real-time 5-minute LMPs from ISO-NE for August 28, 2021 were used as input data to the SM and PM optimization problems \cite{iso_ne_lmp}. Since no reactive power market currently exists, we assumed the Q-LMP to be $10\%$ of the P-LMP \cite{federal_energy_regulatory_commission_payment_2014}. The price ceilings in \cref{eq:opt} were set to be $\overline{\mu}^{iP},\overline{\mu}^{iQ} = 0.2 \; \$/\text{kWh}$, which is almost twice the current average retail rate of 0.129 \$/kWh charged by Eversource, a utility in Massachusetts\footnote{\url{https://www.eversource.com/content/ema-c/residential/my-account/billing-payments/about-your-bill/rates-tariffs/summary-of-electric-rates}}. The overall test feeder was converted to a balanced 3-phase distribution network by (i) assuming switches to be at their normal positions, (ii) converting single phase spot loads to be 3-phase, (iii) assuming cables to be 3-phase transposed, (iv) converting configurations 1 thru 12 to symmetric matrices and (v) modeling shunt capacitors as 3-phase reactive power generators \cite{rabab_tsg}. A PMO was assumed to be at the slack bus, at 13.2kV, with the SMOs at 4.16kV, and each DCA at 120-240V. 

Each SMO was assumed to have anywhere between $\left|\mathcal{N}_{J,i}\right| \in [3,5]$ DCAs with the actual number chosen uniformly at random. The number of DCAs at each SMO $i$ is chosen independently. We set the baseline injections $P^{i0}_j, Q^{i0}_j$ to be equal to the results from the GridLAB-D simulations. Since the injection data was only available up to the primary feeder node level, we artificially disaggregated the injections at each SMO amongst its DCAs, with each DCA being either a net load or net generator. The flexibility bids for the SM $\Delta P_j^i, \Delta Q_j^i$ were also randomly generated, allowing each DCA to offer flexibilities of \textit{up to} $\pm 50\%$ away from their baseline. Thus, the upper and lower limits for the bid flexibilities were set as $\underline{P}^i_j = P^{i0}_j (1 - \underline{\Delta}^j_i), \overline{P}^i_j = P^{i0}_j (1 + \overline{\Delta}^j_i)$, where $\underline{\Delta}^j_i, \overline{\Delta}^j_i \sim \mathcal{U}[0,0.5]$. We focus here on the results for active power only; similar trends were observed for reactive power.

\subsection{SM scheduling \label{sec:sec_mkt}}

\begin{figure}
  \centering
  \begin{subfigure}[b]{\linewidth}
    \includegraphics[width=\linewidth]{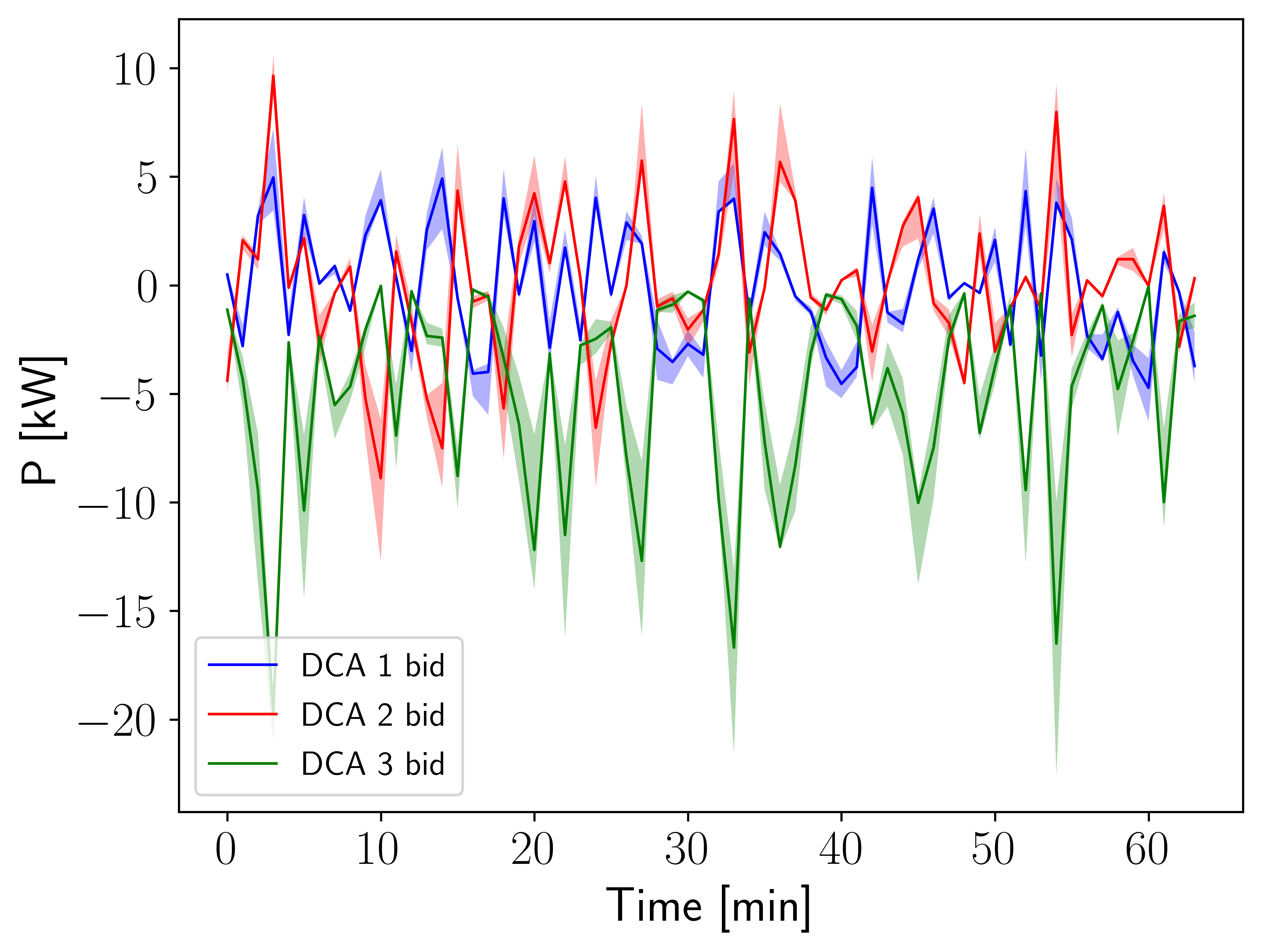}
    \caption{\label{fig:sec_bids} Bids with flexibilities $P^{7,0}_j, \Delta P^7_j$.}
  \end{subfigure}
  \begin{subfigure}[b]{\linewidth}
    \includegraphics[width=\linewidth]{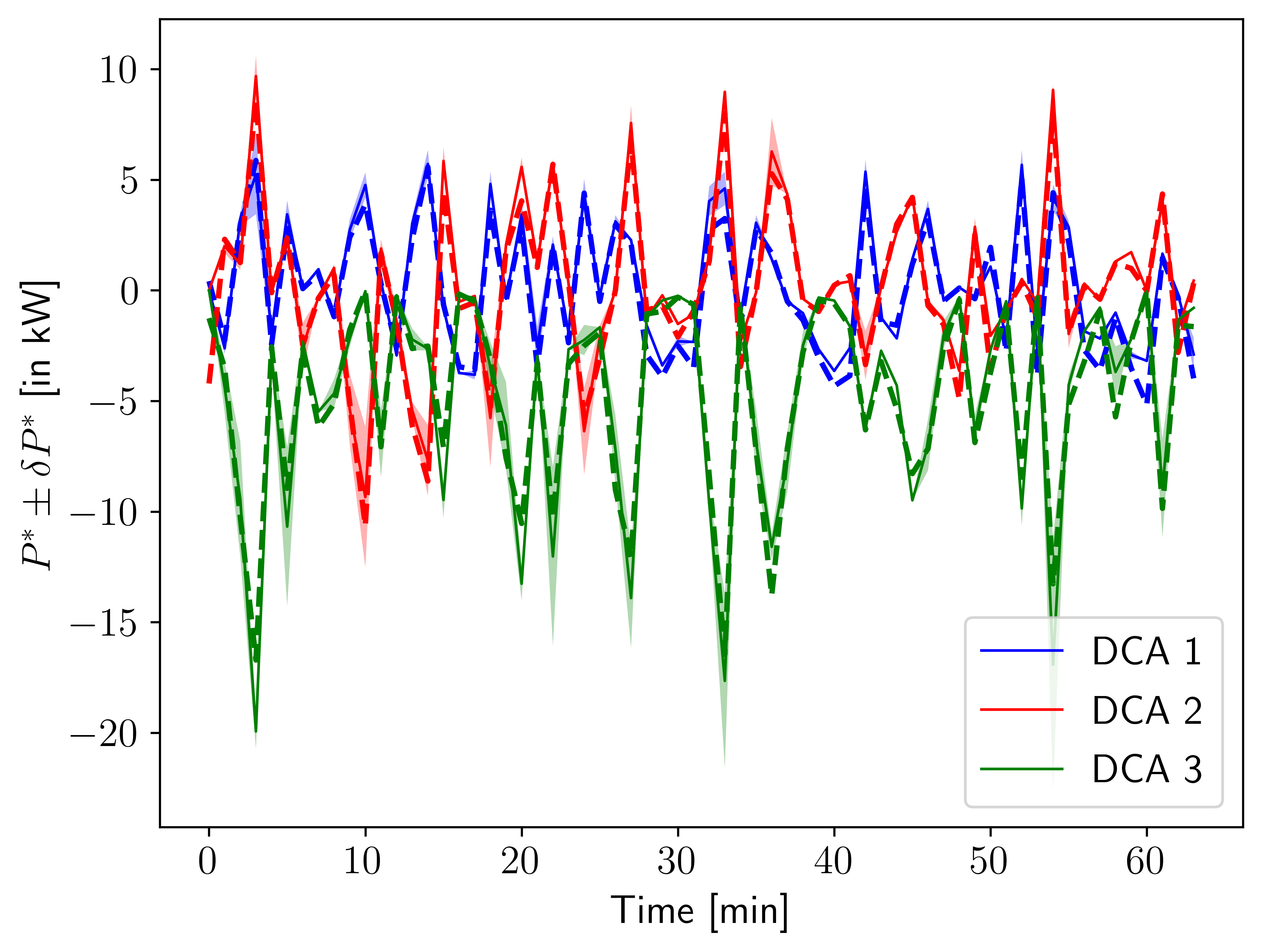}
    \caption{\label{fig:sec_solns} DCA schedules and responses.}
  \end{subfigure}
  \begin{subfigure}[b]{\linewidth}
    \includegraphics[width=\linewidth]{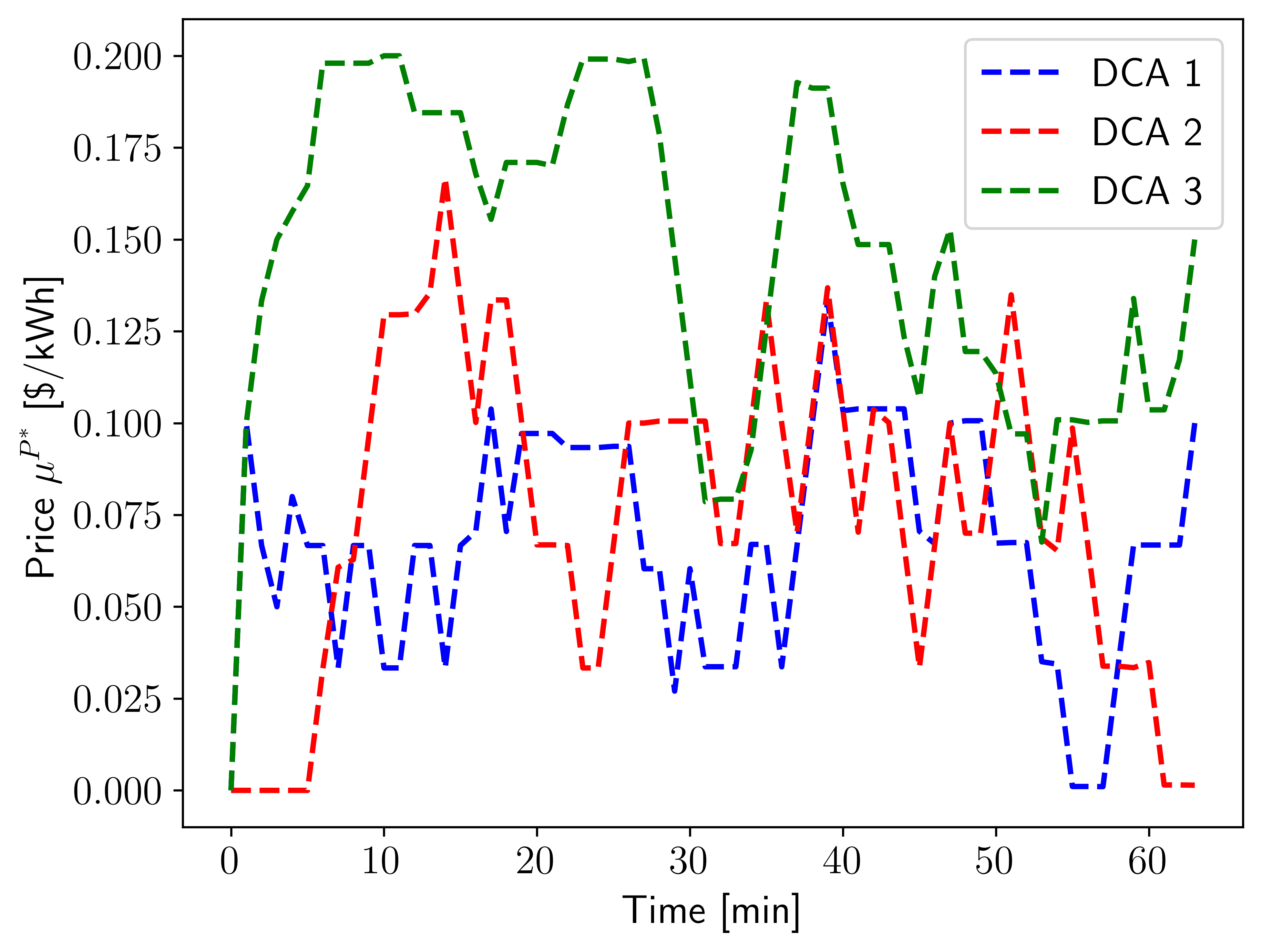}
    \caption{\label{fig:sec_prices} Market cleared local retail tariffs.}
  \end{subfigure}
  \caption{\label{fig:sec_mkt} SM bidding and clearing for primary feeder node 7, with 3 DCAs $j\in \{1,2,3\}$. The solid lines in \cref{fig:sec_bids} and \cref{fig:sec_solns} represent the baseline injection bids and market cleared setpoints, respectively, while the shaded regions around them are the flexibility ranges. Local retail tariffs from the SM $\mu^{7P^*}_j$ are shown in \cref{fig:sec_prices}. The SMO aggregates these PM schedules to bid into the PMO as in \cref{fig:smo_pmo_bids}. The dashed lines in \cref{fig:sec_solns} indicate the actual responses of the DCAs in response to their market cleared schedules.}
\end{figure}

The first step in our use-case study is the SM structure, and its market clearing using the optimization problem outlined in eqs. (\ref{eq:cost})-(\ref{eq:PQbalance}). The bids $\Vec{B}^i_{j}$ corresponding to these parameters are shown in \cref{fig:sec_bids} for a randomly selected SMO $i= 7$ having 3 DCAs $j=1,2,3$. The interval of interest was chosen to be of a 60-min duration, with the actual hour chosen at random. The power injections $P_j^{7^*}$ obtained from solving (\ref{eq:cost})-(\ref{eq:PQbalance}) as well as the corresponding flexibilities, for each DCA $j$, are indicated in \cref{fig:sec_solns}. These two figures clearly illustrate the optimal flexibility range for each of the DCAs, reflecting the ability of the SM to incorporate the constraints of the DCAs, and multiple objectives such as utility, monetary costs, and commitment reliability. The corresponding local electricity tariffs, $\mu_j^{7P^*}$ are shown in \cref{fig:sec_prices} for $j=1,2,3$. Figs. \ref{fig:sec_solns} and \ref{fig:sec_prices} also illustrate the correlation between injections and prices.  For instance, the tariffs for DCA 3 are consistently higher than those for 1 and 2, as DCA 3 is more heavily loaded than the other DCAs. Similarly, tariffs for DCA 1 are lower as its net generation is higher; the price fluctuations are more or less in sync with generation and demand patterns.

\subsection{PM scheduling}

\begin{figure}
  \begin{subfigure}[b]{\linewidth}
    \includegraphics[width=\linewidth]{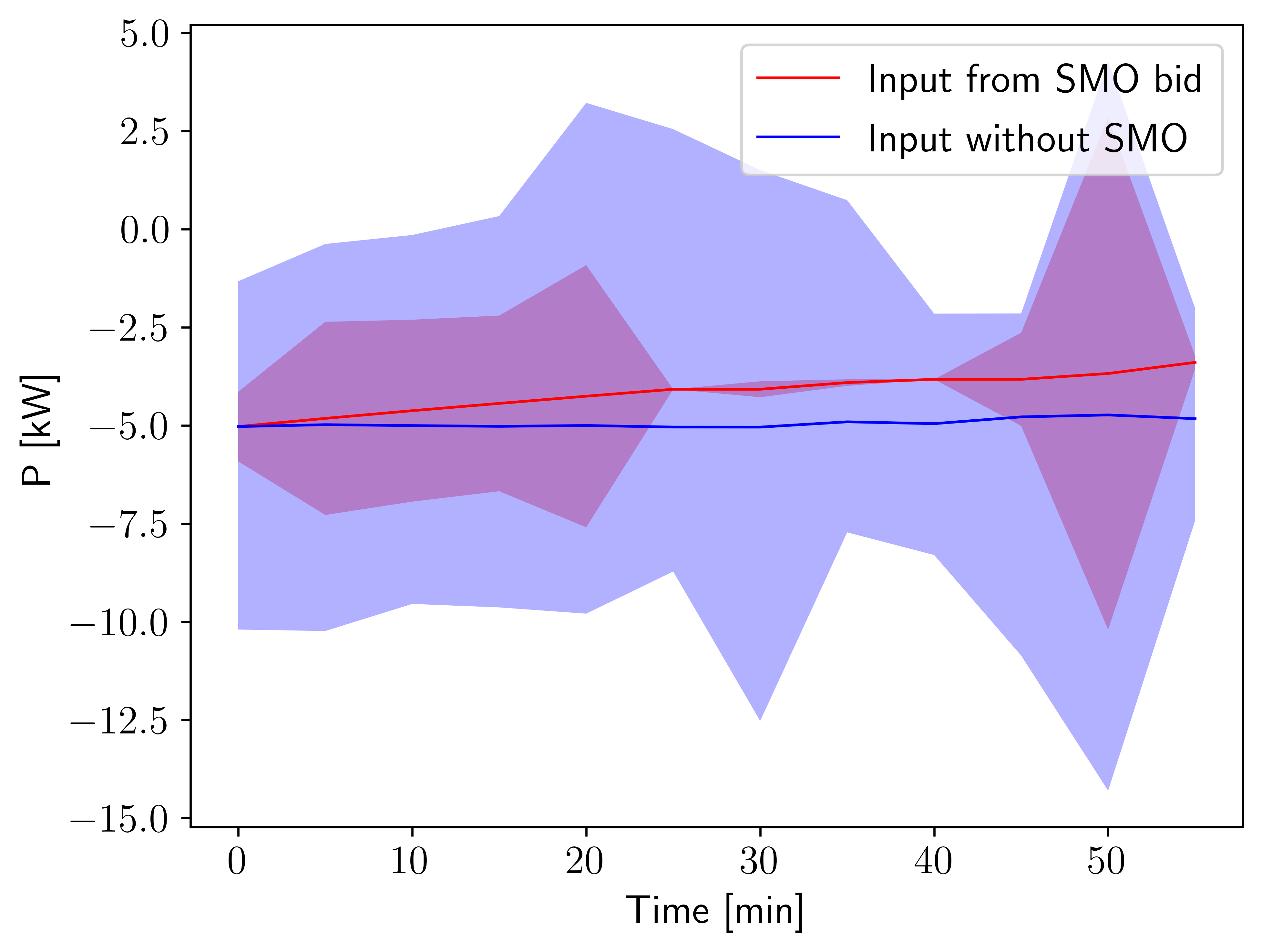}
    \caption{\label{fig:smo_pmo_bids} Inputs from node 7 for PM clearing.}
  \end{subfigure}
   \begin{subfigure}[b]{\linewidth}
    \includegraphics[width=\linewidth]{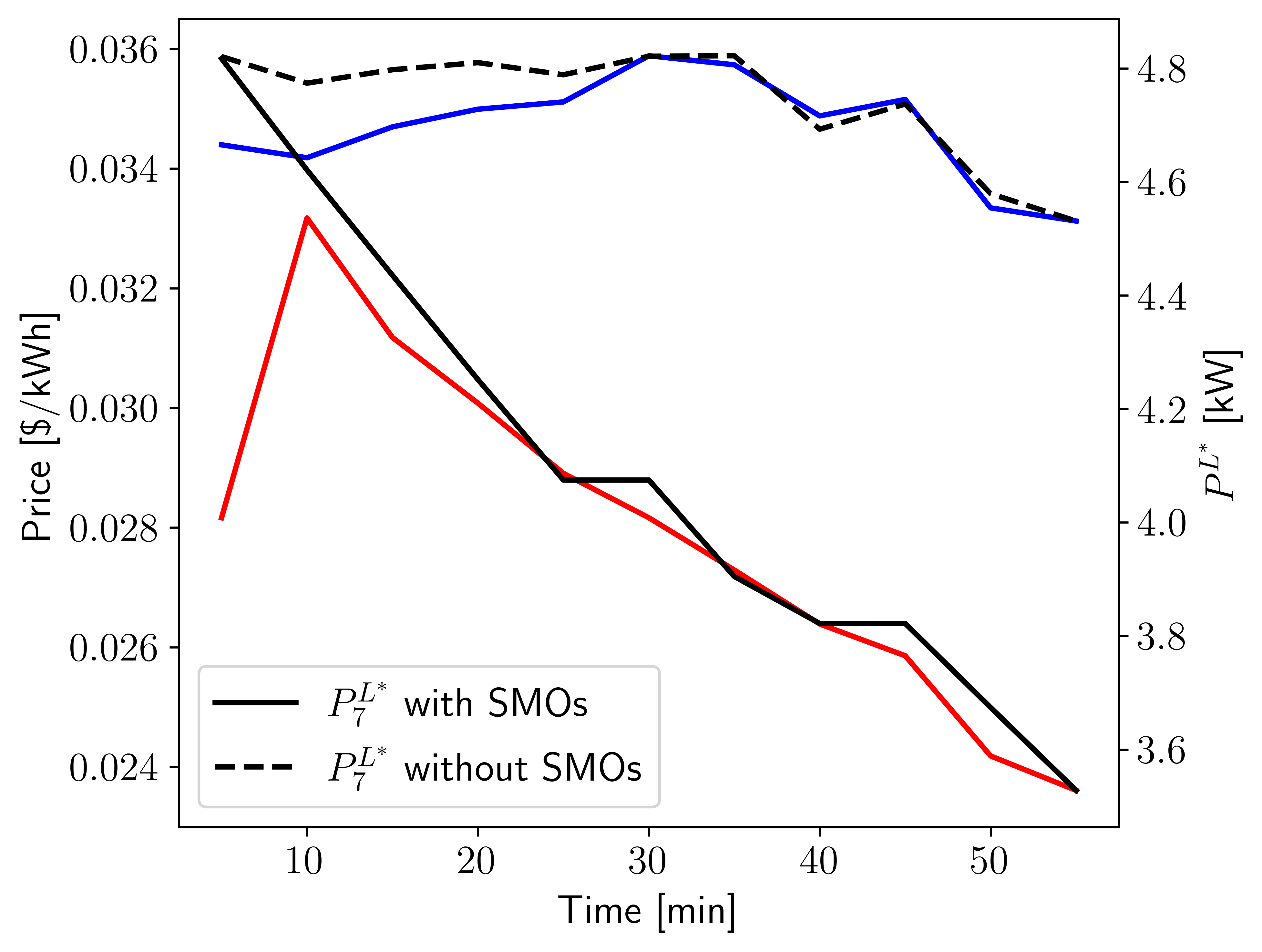}
    \caption{\label{fig:smo_pmo_solns} PM solutions for SMO 7. The solid and dashed black lines are the load injections, while the red and blue lines are the d-LMPs with and without the SMO, respectively.}
  \end{subfigure}
    \begin{subfigure}[b]{\linewidth}
    \includegraphics[width=\linewidth]{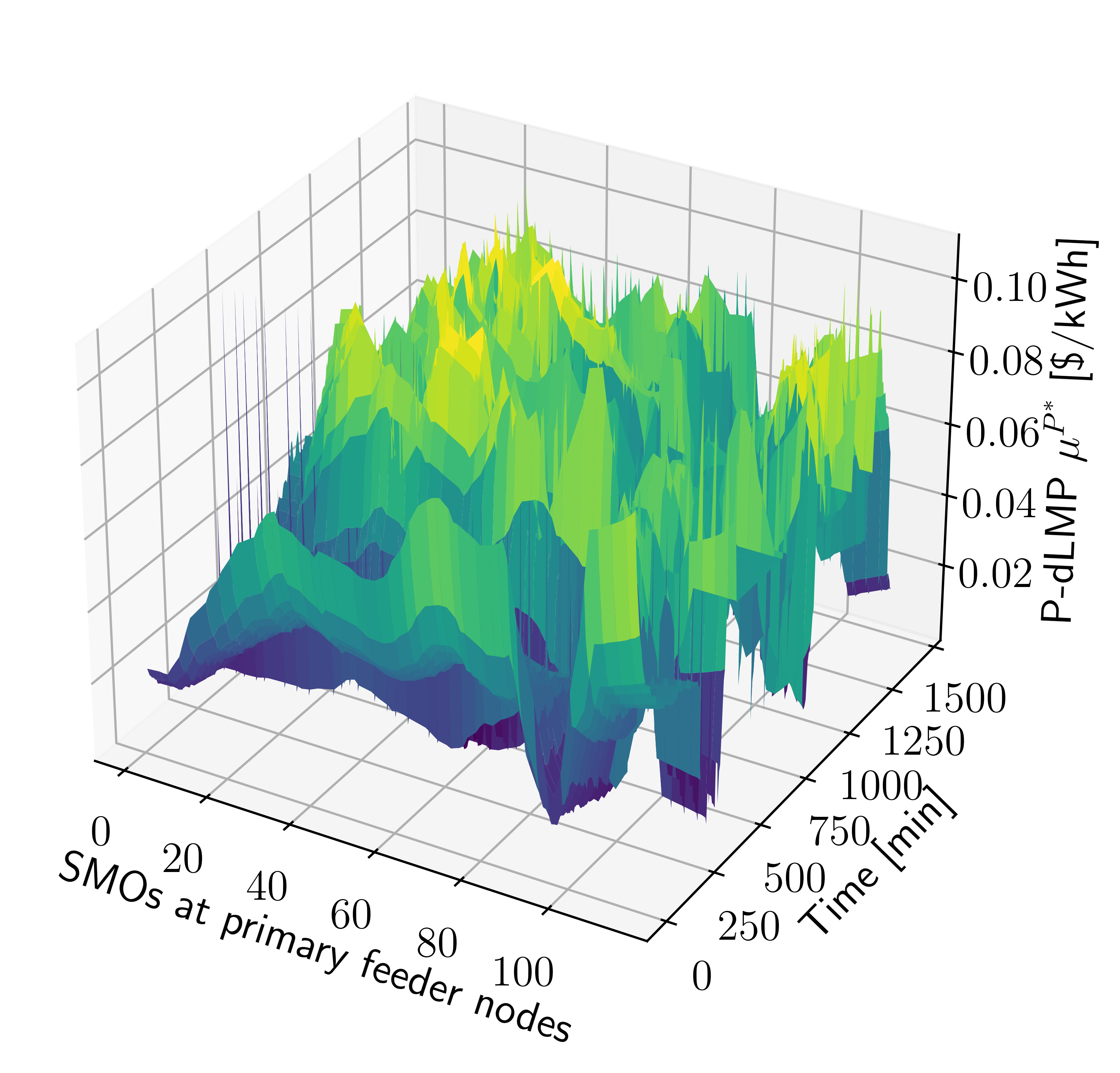}
    \caption{\label{fig:dlmps_surf} d-LMPs across all SMO nodes over 1d.}
  \end{subfigure}
  \caption{\label{fig:pmo_clear} Selected solutions from the PM clearing.}
\end{figure}

The optimal injections with associated flexibilities from the SM clearing in \cref{fig:sec_solns} are aggregated across all three DCAs to form this SMO's bid $P^{0}_7, \Delta P_7$ into the primary level market, as described in \cref{eq:aggregation}. The resulting SMO bids are shown in \cref{fig:smo_pmo_bids}, where the solid red line indicates $P^{0}_7$ and the shaded area indicates the flexibility range $\left[P^{0}_7-\Delta P_7, P^{0}_7+\Delta P_7 \right]$. These bids are in turn used to solve the PM OPF problem in \cref{eq:pmo} using the distributed PAC algorithm, where the SMO's flexible bids $\Delta P_7 = [\underline{P}_7,\overline{P}_7]$ set the feasible operational limits for the power flow constraints in \eqref{eq:pmo}. Solving this optimization problem corresponds to clearing the PM, and determines the PM schedules for the SMO. The results of the PM clearing for SMO $i=7$ are shown in \cref{fig:smo_pmo_solns}. 
Our two-tier market structure generates two sets of schedules and prices, every 1 min and every 5 min for the SM and PM, respectively, shown in \cref{fig:sec_mkt} and \cref{fig:pmo_clear}. We further note that both the local electricity tariffs and the d-LMPs determined by the SM and PM display a high degree of spatio-temporal variations, as shown in \cref{fig:dlmps_surf}. This illustrates the need for local primary and secondary markets to capture such changes with sufficient resolution.

In order to evaluate the impact of the hierarchical structure that we have included in the LEM, we compare the performance of the PM to the case when there is no SM at the lower level. The `without SMO' scenario consists of only a PM, with the PMO directly assuming flexibility ranges for each primary feeder node that best represents an aggregation of local generation and curtailable loads. In what follows, we compare the performance of our hierarchical LEM, i.e., the `with SMO' scenario, with the `without SMO' scenario. First, we compare the inputs into the PM at node 7. 

Figure \ref{fig:smo_pmo_bids} shows that the PMO has a larger flexibility range that may not be accurate or realizable. 
The red curve in \cref{fig:smo_pmo_bids} shows that the flexibility range with SMO is narrower, and reflects the true preferences of the DCAs. Furthermore, the amount of flexibility that the SMO provides to the PMO is also impacted by other factors like the SM retail costs and the commitment scores of each of its DCAs, both of which vary with time. As a result, the `with SMO' case is more performant as the baseline injection is optimized in comparison to the relatively ad-hoc choice in the without SMO case (the blue curve in \cref{fig:smo_pmo_bids}). 

\begin{figure}[!htbp]
  \centering
  \begin{subfigure}[b]{\linewidth}
    \includegraphics[width=\linewidth]{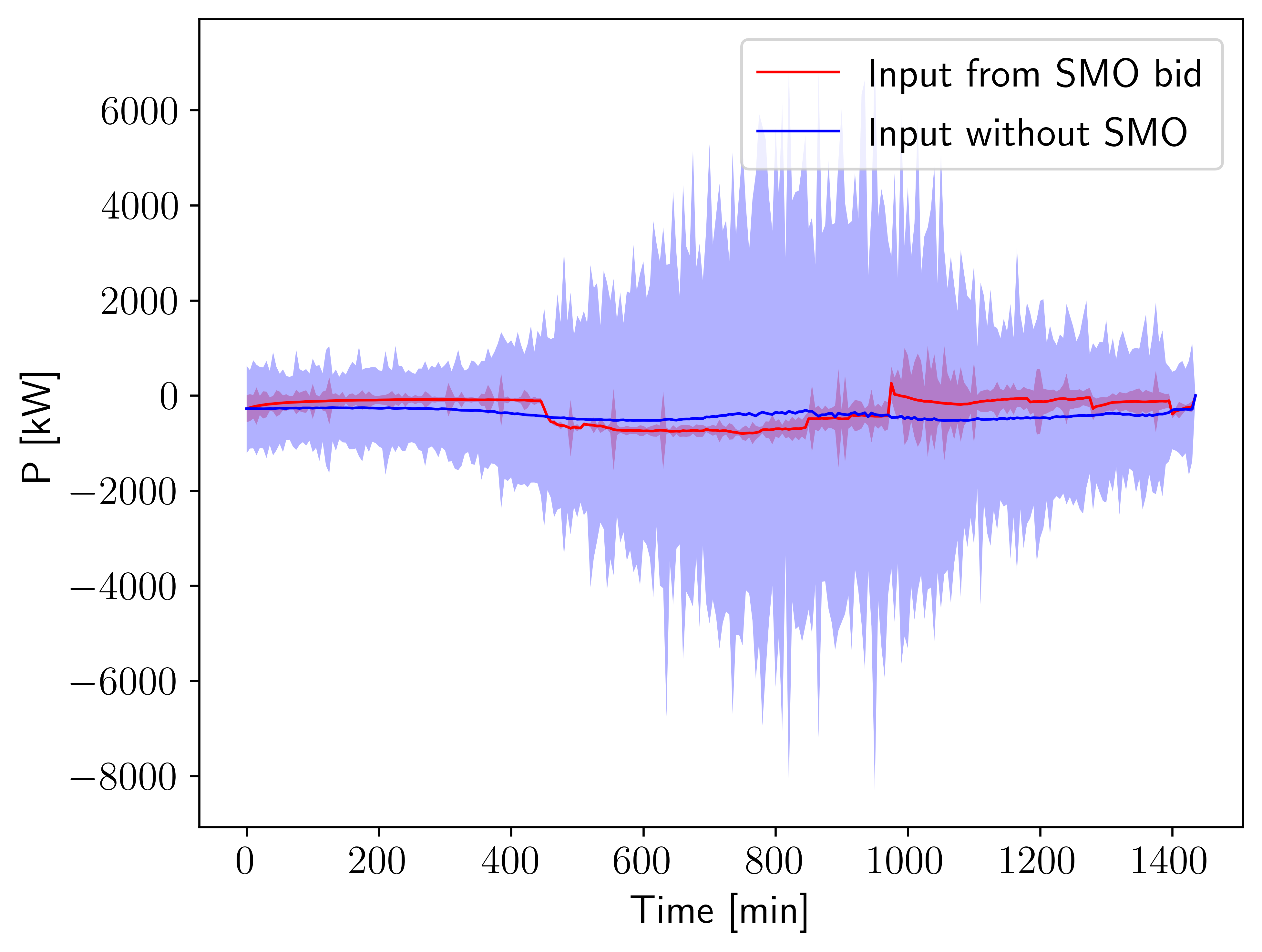}
    \caption{\label{fig:smo_pmo_bids_all} Inputs to PM aggregated across all primary feeder nodes except the slack bus.}
  \end{subfigure}
  \begin{subfigure}[b]{\linewidth}
    \includegraphics[width=\linewidth]{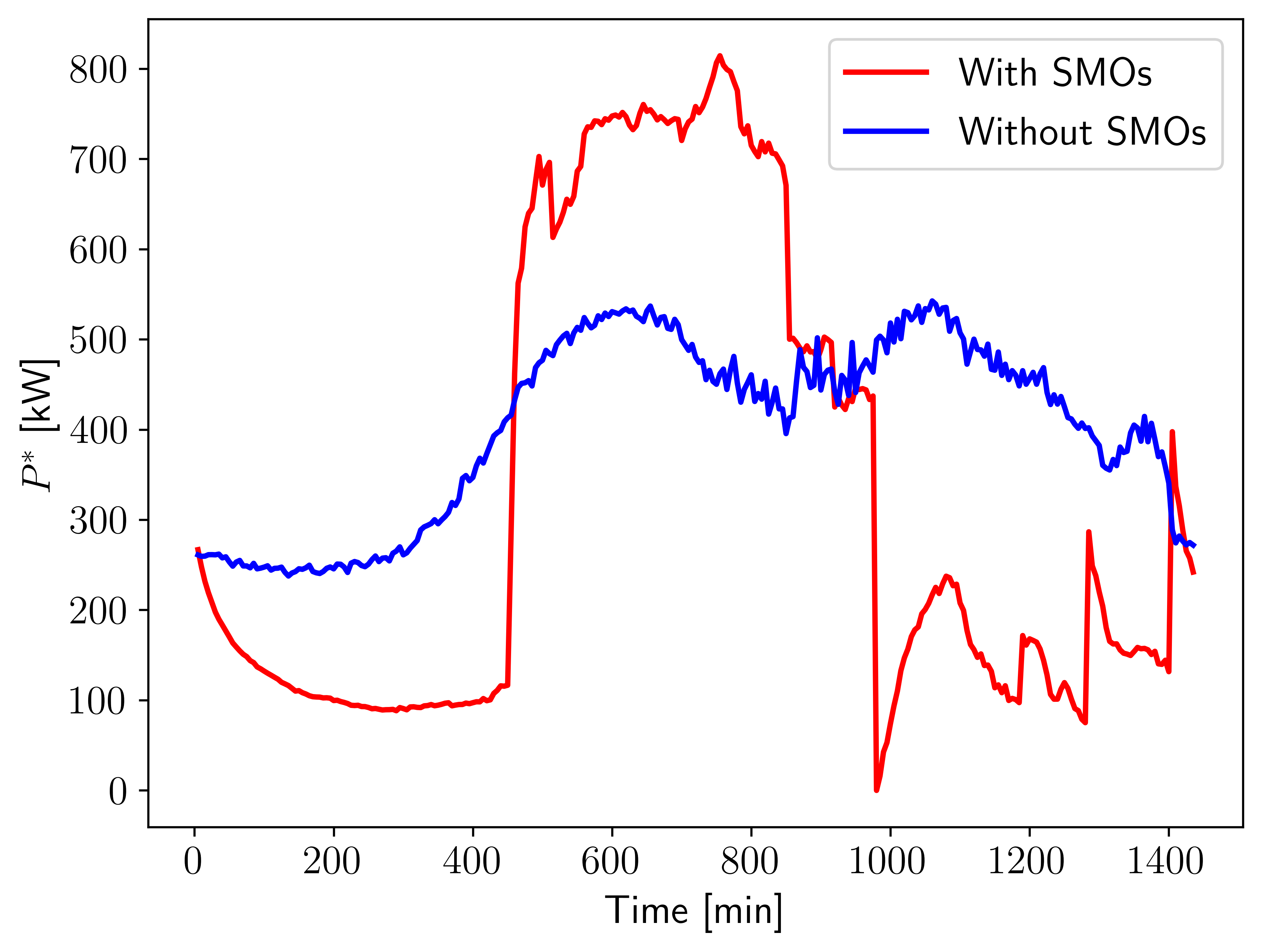}
    \caption{\label{fig:smo_pmo_Pslack} PM solutions for net injections at the slack bus.}
  \end{subfigure}
  \caption{Comparison of PM bids (or inputs) and slack bus injections, with and without SM. The slack bus (node 149) is connected to the substation and distribution transformer. Positive injections here indicate that the feeder as a whole is importing power from the main grid.}
\end{figure}

\begin{figure}[!htbp]
  \centering
  \begin{subfigure}[b]{\linewidth}
    \includegraphics[width=\linewidth]{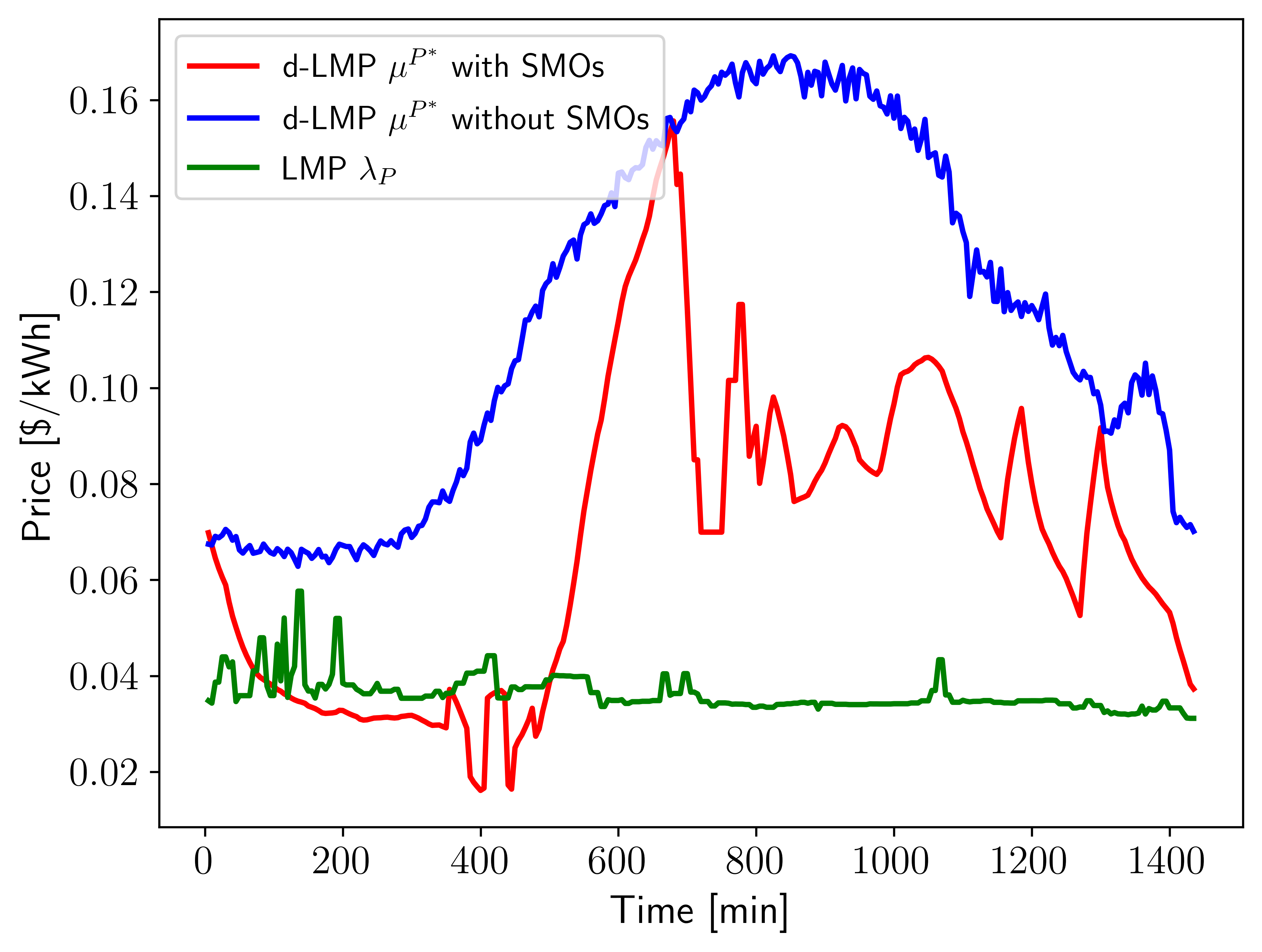}
    \caption{\label{fig:smo_pmo_dLMPs} d-LMPs averaged across all 79 primary nodes considered, with and without the SM, compared to the LMP from the WEM.}
  \end{subfigure}
  \begin{subfigure}[b]{\linewidth}
    \includegraphics[width=\linewidth]{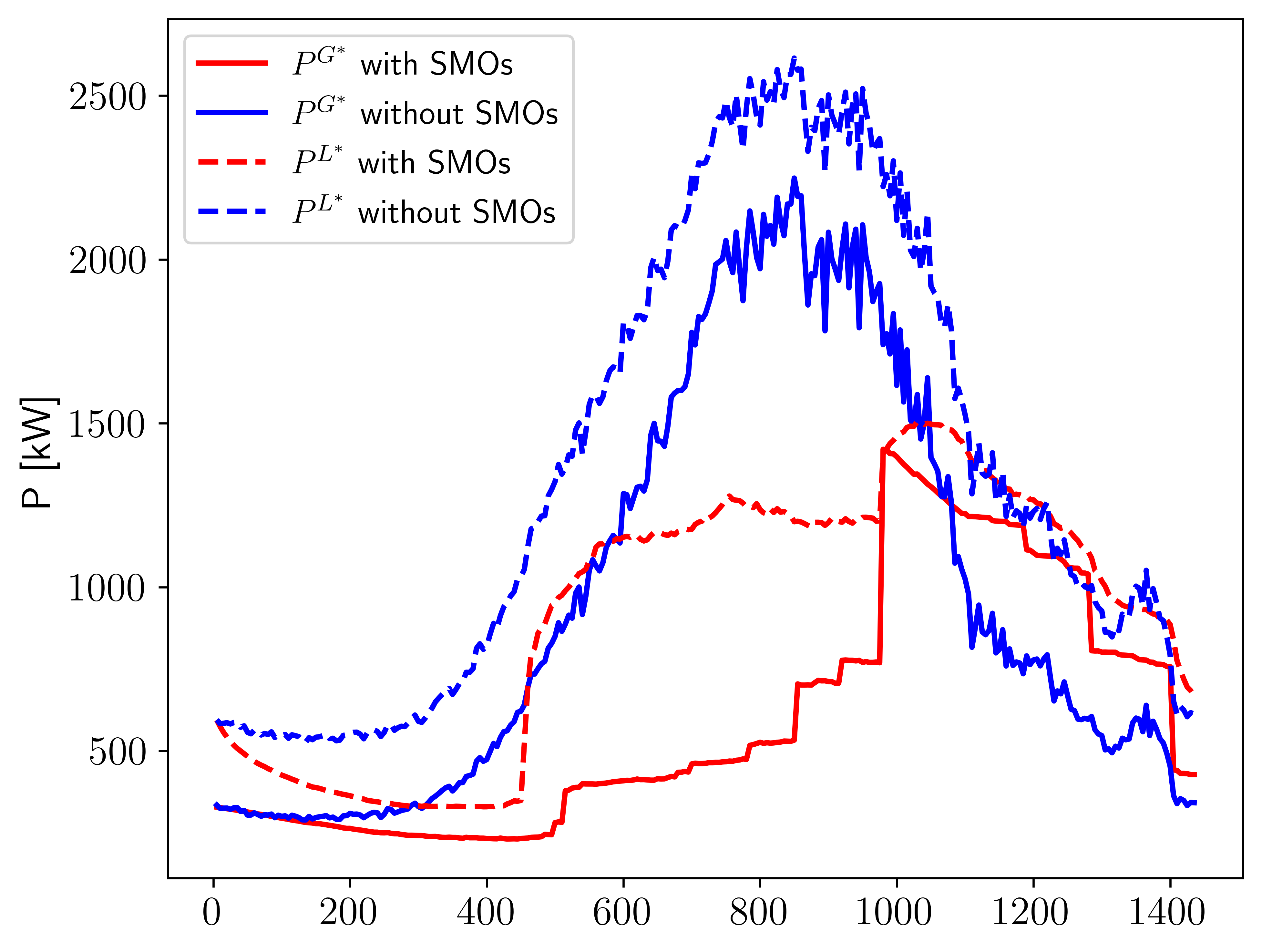}
    \caption{\label{fig:smo_pmo_PL_PG} PM load and generation injections, summed over all primary feeder nodes except the slack bus.}
  \end{subfigure}
  \caption{Comparison of PM solutions obtained with and without SM.}
\end{figure}

We next compare the performance of our hierarchical market across the entire primary feeder consisting of all 79 SMO nodes, over the course of the whole simulation period of 24 hours. In \cref{fig:smo_pmo_bids_all}, the inputs to the PMO are shown (the red curve), with all SMO solutions aggregated across all 79 primary feeder nodes $i \in \mathcal{N}_I$ and for the entire day. We see that without the additional visibility and granularity offered by the SM structure, the PM would assume much larger ranges for the injection limits in the `without SMO' case (the blue curve) when compared to the `with SMO' case. These are less accurate and may also be overoptimistic in terms of how much flexibility can be realistically expected from the DCAs, which in turn can cause issues in case of reneged commitments. It should be pointed out that the amounts of local generation seen in \cref{fig:pmo_clear} and \cref{fig:smo_pmo_PL_PG} are above the installed PV capacity of 510.3 kW. This is because while generating the synthetic flexibility bids for the DCAs, we allowed for the possibility of additional DERs like batteries, EVs and curtailable or shiftable loads, present at each of these secondary feeders, which weren't explicitly modeled in the GridLAB-D simulation.

In \cref{fig:smo_pmo_dLMPs}, we see that the d-LMPs both with and without the SMO are generally higher than the LMP, which is expected since the d-LMPs account for additional costs associated with congestion, line losses and other delivery charges incurred by the PMO and DSO in the distribution network, downstream of the substation. The d-LMP with the SMO does fall slightly below the LMP between 100-500 minutes (02:00:0700). This can be explained by the total electricity demand being low during this period which in turn occurs as the SMOs are able to curtail flexible loads to a larger extent by coordinating their DCAs more intelligently and compensate them accordingly at the local retail tariff rate. In fact, we find that the SMOs are able to achieve higher levels of load curtailment throughout the course of the day when compared to the case without SMOs. Once again, this is likely because the SMO can access additional information on DCA's preferences and effectively utilize any additional flexibility that they're willing to provide. The SM allows the SMO to more efficiently allocate resources amongst the secondary feeders at each primary feeder node, and take advantage of differences in load and generation profiles across DCAs over time since they could potentially complement each other.  


The second observation from  \cref{fig:smo_pmo_dLMPs} and \cref{fig:smo_pmo_PL_PG} is that the `with SMO' case schedules lower levels of local generation mid-day compared to the `without SMO' case. This may be due to a combination of multiple objectives utilized in the SM that include both net costs and flexibility. The optimal behavior as a result, as predicted by the LEM, is one where more power is purchased from the main transmission grid rather than from local generation mid-day. This is also supported by \cref{fig:smo_pmo_dLMPs} which shows that such a behavior leads to lower d-LMPs and reduced distribution network costs with the hierarchical LEM than without the SMO. This is desirable since the SMOs can then reduce the retail tariff charged to their DCAs, improving affordability for customers, as seen in \cref{tab:summary}. It also ensures that DSOs aren't over-compensating prosumers with DERs. This can help avoid excessive cross-subsidies from consumers to prosumers which is a major challenge associated with net energy metering (NEM) programs today \cite{Eid2014TheObjectives,Picciariello2015ElectricityProsumers}, and can thus produce more equitable allocations.

Figures \ref{fig:smo_pmo_Pslack}, \ref{fig:smo_pmo_dLMPs}, and \ref{fig:smo_pmo_PL_PG}  correspond to the main conclusions of the proposed LEM. In all three figures, the red curves correspond to the behavior with the LEM while the blue curves correspond to the `without SMO' case. The red curve in \cref{fig:smo_pmo_Pslack} shows that the LEM schedules generation from the bulk grid more in the middle of the day and less otherwise; those in \cref{fig:smo_pmo_PL_PG} show that it's advantageous to increase local generation in the latter part of the day and to curtail load in the earlier part of the day. The LEM determines that the IEEE-123 feeder needs to import around 700 kW between minute 400 to minute 850, and less than 300 kW from minute 1000 onward. This behavior is significantly different from the market structure without SMOs, as the primary market alone does not have the granular customer level information to accurately estimate the power injections and their associated flexibilities. Finally, \cref{fig:smo_pmo_dLMPs} shows the optimal d-LMP from the LEM that enables the overall generation mix as shown in \cref{fig:smo_pmo_Pslack} and \cref{fig:smo_pmo_PL_PG}, and that it is lower than what the `without SMO' case predicts. 

\begin{table}
\centering
\caption{\label{tab:summary} Summary financial metrics for our simulations under different types of market structures.}
\begin{tabular}{@{}lccc@{}}
\toprule
                                         & \textbf{SM + PM} & \textbf{PM only} & \textbf{No LEM} \\ \midrule
\textbf{Avg. P d-LMP {[}\$/kWh{]}}        & 0.064                     & 0.116                     & N/A             \\
\textbf{Avg retail tariff {[}\$/kWh{]}}  & 0.082                     & 0.116                     & 0.129           \\
\bottomrule
\end{tabular}
\end{table}


\section{Concluding Remarks \label{sec:conclusion}}
A hierarchical local electricity market (LEM) structure was proposed in this paper with a secondary market (SM) at the lower level representing secondary feeders and a primary market (PM) at the upper level, representing primary feeders, in order to effectively use DERs. The lower level SM enforces budget, power balance and flexibility constraints and accounts for costs related to consumers, such as their disutility, flexibility limits, and commitment reliability, while the upper level PM enforces grid physics constraints such as power balance and capacity limits, and also minimizes line losses. The hierarchical LEM is evaluated using a modified IEEE-123 bus with high DER penetration, with each primary feeder consisting of several secondary feeders. Realistic power injections and load profiles were obtained over the course of 24 hours from GridLAB-D. The performance of the LEM was illustrated by delineating the family of power-injection profiles across the primary and secondary feeders as well as the corresponding local electricity tariffs that vary across the distribution grid. It was shown that the overall LEM is capable of capturing fine-grain variations across the primary feeders and even further across secondary feeders so that the power injections and corresponding variable tariffs accurately charge or compensate DERs, capture consumer flexibilities, DER capabilities and constraints, as well as constraints and costs stemming from power physics. 

The proposed hierarchical LEM represents the first step in formulating a market structure that allows disparate DER assets to participate and be appropriately compensated. Several other steps needs to be executed to develop a complete retail market with various products. First, multiperiod extensions of our optimization frameworks at both the SM and PM levels need to be carried out. Next, the co-optimization of our SM and PM in the real-time market needs to be addressed, along with other markets for ancillary services and reserves. In doing so, we also hope to fully address the issues around settlement and billing in our hierarchical LEM, as well as the handling of unmet or reneged commitments in real-time. Also required is the development of advanced game theoretic approaches that could be used by the SMO to generate its bids into the PM from the SM solutions, instead of the simple direct sum aggregation used at present. Similarly, the possibility of strategic bidding by DCAs in the SM, and methods to counter this needs to be examined as well. This will also help guide the design of a consumer-level market within each secondary feeder, forming the final tier of the proposed hierarchical LEM. Finally, more realistic distribution-level test cases and datasets are planned to be developed to validate the overall LEM, leveraging both simulations as well as real-world data \cite{postigo2017review,bu_time-series_2019}.

\section*{Acknowledgments}
We gratefully acknowledge several useful discussions with Prof. Anurag Srivastava, Dr. Amit Chakraborty, and Dr. Biswadip Dey.

\bibliographystyle{IEEEtran}
\bibliography{refs, mendeley_vineet, zotero_vineet, references_mendeley}

\end{document}